\RequirePackage{fix-cm}
\documentclass[smallcondensed]{svjour3}     %
\usepackage{xcolor}

\usepackage{graphicx}

\usepackage{amssymb}
\usepackage{amsmath}
\usepackage{algorithm}
\usepackage{algorithmic}
\usepackage{mathrsfs}
\usepackage{cite}
\usepackage{mathtools}

\usepackage{fixme}
\fxsetup{layout=inline, theme=colorsig}

\FXRegisterAuthor{gs}{ags}{GS}
\FXRegisterAuthor{lcw}{alcw}{LCW}

\usepackage[squaren]{SIunits}

\usepackage{lineno}

\newcommand*\patchAmsMathEnvironmentForLineno[1]{%
  \expandafter\let\csname old#1\expandafter\endcsname\csname #1\endcsname
  \expandafter\let\csname oldend#1\expandafter\endcsname\csname end#1\endcsname
  \renewenvironment{#1}%
     {\linenomath\csname old#1\endcsname}%
     {\csname oldend#1\endcsname\endlinenomath}}%
\newcommand*\patchBothAmsMathEnvironmentsForLineno[1]{%
  \patchAmsMathEnvironmentForLineno{#1}%
  \patchAmsMathEnvironmentForLineno{#1*}}%
\AtBeginDocument{%
\patchBothAmsMathEnvironmentsForLineno{equation}%
\patchBothAmsMathEnvironmentsForLineno{align}%
\patchBothAmsMathEnvironmentsForLineno{flalign}%
\patchBothAmsMathEnvironmentsForLineno{alignat}%
\patchBothAmsMathEnvironmentsForLineno{gather}%
\patchBothAmsMathEnvironmentsForLineno{multline}%
}

\usepackage[colorlinks,citecolor=blue!80!black,urlcolor=blue,linkcolor=red!80!black]{hyperref}
\usepackage{tikz}
\usetikzlibrary{matrix,arrows}
\usepackage{pgfplots}

\DeclareMathOperator{\sym}{sym}
\DeclareMathOperator{\trace}{tr}
\DeclareMathOperator{\diag}{diag}

\newcommand{\bs}[1]{\boldsymbol{#1}}

\newcommand{\point}[1] {\ensuremath{\boldsymbol{#1}}}
\renewcommand{\vec}[1] {\ensuremath{\boldsymbol{#1}}}

\newcommand{\ten}[1] {\ensuremath{\boldsymbol{#1}}}
\newcommand{\tenfour}[1] {\ensuremath{\boldsymbol{\mathsf{#1}}}}

\newcommand{\vv}{\ensuremath{\vec{v}}}
\newcommand{\xx}{\ensuremath{\point{x}}}
\newcommand{\ff}{\ensuremath{\vec{f}}}
\newcommand{\FF}{\ensuremath{\mathbf{F}}}
\renewcommand{\gg}{\ensuremath{\vec{g}}}
\newcommand{\qq}{\ensuremath{\vec{q}}}
\newcommand{\pp}{\ensuremath{\vec{p}}}
\newcommand{\tvv}{\ensuremath{\tilde{\vec{v}}}}
\newcommand{\ww}{\ensuremath{\vec{w}}}

\newcommand{\HH}{\ensuremath{\ten{H}}}
\newcommand{\EE}{\ensuremath{\ten{E}}}
\newcommand{\tEE}{\ensuremath{\tilde{\ten{E}}}}
\newcommand{\nn}{\ensuremath{\vec{n}}}
\renewcommand{\SS}{\ensuremath{\ten{S}}}

\newcommand{\cc}{\ensuremath{\mathbf{c}}}
\newcommand{\cp}{\ensuremath{{c_p}}}
\newcommand{\cs}{\ensuremath{{c_s}}}
\newcommand{\tcp}{\ensuremath{{\tilde{c}_p}}}

\newcommand{\nxnx}[1] {\ensuremath{\nn^-\!\times\!\left(\nn^-\!\times{#1}\right)}}

\newcommand{\Grad} {\ensuremath{\nabla}}
\newcommand{\Div} {\ensuremath{\nabla\cdot}}
\newcommand{\Curl} {\ensuremath{\nabla\times}}
\newcommand{\CC}{\ensuremath{\tenfour{C}}}

\newcommand{\mXY}{\ensuremath{P^h\times Q^h}}
\newcommand{\mEE}{\ensuremath{\vec{E}}}
\newcommand{\mHH}{\ensuremath{\vec{H}}}
\newcommand{\mFF}{\ensuremath{\vec{F}}}
\newcommand{\mGG}{\ensuremath{\vec{G}}}

\newcommand{\mJs}{\ensuremath{\vec{J}_{s}}}
\newcommand{\mg}{\ensuremath{\vec{g}}}

\newcommand{\aeXY}{\ensuremath{P^h\times Q^h}}

\newcommand{\eval}[2][\right]{\relax
  \ifx#1\right\relax \left.\fi#2#1\rvert}

\DeclarePairedDelimiter\abs{\lvert}{\rvert}%
\DeclarePairedDelimiter\norm{\lVert}{\rVert}%

\makeatletter
\let\oldabs\abs
\def\abs{\@ifstar{\oldabs}{\oldabs*}}
\let\oldnorm\norm
\def\norm{\@ifstar{\oldnorm}{\oldnorm*}}
\makeatother

\newcommand{\jump}[1] {\ensuremath{[\![#1]\!]}}
\newcommand{\diff}[1] {\ensuremath{[#1]}}
\newcommand{\mean}[1] {\ensuremath{\{\!\!\{#1\}\!\!\}}}

\newcommand{\De} {\ensuremath{{\Omega^e}}}
\newcommand{\Dep} {\ensuremath{{\Omega^{e'}}}}

\newcommand{\sR}{\mathbb{R}}

\begin{document}
\title{Discretely exact derivatives for hyperbolic PDE-constrained
  optimization problems discretized by the discontinuous Galerkin
  method\footnote{This document has been approved for public release;
    its distribution is unlimited.}}

\subtitle{}
\titlerunning{Derivatives for dG-discretized hyperbolic optimization problems}
\authorrunning{L.\ C.\ Wilcox, G.\ Stadler, T.\ Bui-Thanh and O.\ Ghattas}

\author{Lucas C.~Wilcox$^1$, Georg Stadler$^2$,\\ Tan Bui-Thanh$^{2,3}$ and Omar Ghattas$^{2,4}$}

\institute{$^1${Department of Applied Mathematics, Naval Postgraduate School, Monterey, CA}\\
$^2$Institute for Computational Engineering \& Sciences, The
  University of Texas at Austin, Austin, TX\\
$^3$Department of Aerospace Engineering \& Engineering Mechanics, The
  University of Texas at Austin, Austin, TX\\
$^4$Departments of Mechanical Engineering and Jackson School of Geosciences, The
  University of Texas at Austin, Austin, TX
}

\maketitle
\begin{abstract}
This paper discusses the computation of derivatives for optimization
problems governed by linear hyperbolic systems of partial differential
equations (PDEs) that are discretized by the discontinuous Galerkin
(dG) method. An efficient and accurate computation of these
derivatives is important, for instance, in inverse problems
and optimal control problems.
This computation is usually based on an adjoint PDE system, and the
question addressed in this paper is how the discretization of this
adjoint system should relate to the dG discretization of the
hyperbolic state equation. Adjoint-based derivatives can either be computed
before or after discretization; these two options are often referred
to as the optimize-then-discretize and discretize-then-optimize
approaches. We discuss the relation between these two options for dG
discretizations in space and Runge--Kutta time integration.  The
influence of different dG formulations and of numerical quadrature is
discussed. Discretely exact discretizations for several hyperbolic
optimization problems are derived, including the advection equation,
Maxwell's equations and the coupled elastic-acoustic wave equation.
We find that the discrete adjoint equation inherits a natural dG
discretization from the discretization of the state equation and
that the expressions for the discretely exact gradient often have to
take into account contributions from element faces.  For the coupled
elastic-acoustic wave equation, the correctness and accuracy of our
derivative expressions are illustrated by comparisons with finite
difference gradients. The results show that a straightforward
discretization of the continuous gradient differs from the discretely
exact gradient, and thus is not consistent with the discretized
objective. This inconsistency may cause difficulties in the
convergence of gradient based algorithms for solving optimization
problems.

\end{abstract}
\keywords{%
  Discontinuous Galerkin \and
  PDE-constrained optimization \and
  Discrete adjoints \and
  Elastic wave equation \and
  Maxwell's equations}

\section{Introduction}

Derivatives of functionals, whose evaluation depends on the solution
of a partial differential equation (PDE), are required, for instance,
in inverse problems and optimal control problems, and play a role in
error analysis and a posteriori error estimation.  An efficient method
to compute derivatives of functionals that require the solution of a
state PDE is through the solution of an adjoint equation. In
general, this adjoint PDE differs from the state PDE\@.
For instance, for a state equation that involves a first-order time
derivative, the adjoint equation must be solved backwards in time.

If the partial differential equation is hyperbolic, the discontinuous
Galerkin (dG) method is often a good choice to approximate the
solution due to its stability properties, flexibility, accuracy and
ease of parallelization.  Having chosen a dG discretization for the state PDE,
the question arises how to discretize the adjoint equation and the
expression for the gradient, and whether and how their discretization should be
related to the discretization of the state equation.  One approach
is to discretize the adjoint equation and the gradient independently
from the state equation, possibly leading to inaccurate derivatives
as discussed below. A different approach is to derive the discrete adjoint
equation based on the discretized state PDE and a discretization of
the cost functional. Sometimes this latter approach is called the
discretize-then-optimize approach, while the former is known as
optimize-then-discretize. For standard Galerkin discretizations, these
two possibilities usually coincide; however, they can differ, for instance, for
stabilized finite element methods and for shape
derivatives~\cite{CollisHeinkenschloss02, Gunzburger03,
  HinzePinnauUlbrichEtAl09, Braack09}.  In this paper, we study the interplay
between these issues for derivative computation in optimization
problems and discretization by the discontinuous Galerkin method.

While computing derivatives through adjoints on the
infinite-dimensional level and then discretizing the resulting
expressions (optimize-then-discretize) seems convenient, this approach
can lead to inaccurate gradients that are not proper derivatives of
any optimization problem.
This can lead to convergence problems in
optimization algorithms due to inconsistencies between the cost
functional and gradients~\cite{Gunzburger03}. This inaccuracy is
amplified when inconsistent gradients are used to approximate second
derivatives based on first derivatives, as in quasi-Newton methods such as the BFGS method.
Discretizing the PDE and the cost functional first
(discretize-then-optimize), and then computing the (discrete)
derivatives guarantees consistency. However, when an advanced
discretization method is used, computing the discrete derivatives can
be challenging.  Thus, understanding the relation between
discretization and adjoint-based derivative computation is
important. In this paper, we compute derivatives based on the
discretized equation and then study how the resulting adjoint
discretization relates to the dG discretization of the state
equation, and study the corresponding consistency issues for the
gradient.

\emph{Related work:} Discretely exact gradients can also be generated
via algorithmic differentiation (AD)~\cite{GriewankWalther08}. While
AD guarantees the computation of exact discrete gradients, it is
usually slower than hand-coded derivatives. Moreover, applying AD to
parallel implementations can be challenging~\cite{UtkeHascoetHeimbachEtAl09}.
The notion of \emph{adjoint consistency} for dG (see~\cite{Hartmann07, AlexeSandu10,
  HarrimanGavaghanSuli04, OliverDarmofal09, SchutzMay13}) is related to the discussion in this
paper. Adjoint consistency refers to the fact that the exact solution
of the dual (or adjoint) problem satisfies the discrete adjoint
equation. This property is important for dG discretizations to obtain
optimal-order $L^2$-convergence with respect to target
functionals. The focus of this paper goes beyond adjoint
consistency to consider 
consistency of the gradient expressions, and considers in what sense
the discrete gradient is a discretization of the continuous
gradient. For discontinuous Galerkin discretization, the latter aspect
is called \textit{dual consistency} in~\cite{AlexeSandu10}. A
systematic study presented in~\cite{Leykekhman12} compares different
dG methods for linear-quadratic optimal control problems subject to
advection-diffusion-reaction equations. In particular, the author
targets commutative dG schemes, i.e., schemes for which
dG discretization and the gradient derivation commute. Error estimates
and numerical experiments illustrate that commutative schemes have
desirable properties for optimal control problems.

\emph{Contributions:} Using example problems, we illustrate that the
discrete adjoint of a dG discretization is, again, a dG discretization
of the continuous adjoint equation. In particular, an upwind numerical
flux for the hyperbolic state equation turns into a downwind flux in
the adjoint, which has to be solved backwards in time and converges at
the same convergence order as the state equation. We discuss the
implications of numerical quadrature and of the choice of the weak or
strong form of the dG discretization on the adjoint system.  In our
examples, we illustrate the computation of derivatives with respect to
parameter fields entering in the hyperbolic system either as a
coefficient or as forcing terms. Moreover, we show that discretely
exact gradients often involve contributions at element faces, which
are likely to be neglected in an optimize-then-discretize
approach. These contributions are a consequence of the discontinuous
basis functions employed in the dG method and since they are at the
order of the discretization error, they are particularly important for
not fully resolved problems.

\emph{Limitations:} We restrict ourselves to problems
governed by \emph{linear} hyperbolic systems. This allows for an explicit
computation of the upwind numerical flux in the dG method through the
solution of a Riemann problem. Linear problems usually do not
require flux limiting and do not develop shocks in the solution, which 
makes  the computation of derivatives problematic
since numerical fluxes with limiters are often non-differentiable and
defining adjoints when the state solution involves shocks is
a challenge~\cite{GilesUlbrich10,GilesUlbrich10a}.

\emph{Organization:} Next, in Section~\ref{sec:preliminaries}, we
discuss the interplay of the derivative computation of a cost
functional with the spatial and temporal discretization of the
governing hyperbolic system; moreover, we discuss the effects of
numerical quadrature. For examples of linear hyperbolic systems with
increasing complexity we derive the discrete adjoint systems and
gradients in Section~\ref{sec:examples}, and we summarize important
observations. In Section~\ref{sec:numerics}, we numerically verify our
expressions for the discretely exact gradient for a cost functional
involving the coupled acoustic-elastic wave equation by comparing to
finite differences, and finally, in Section~\ref{sec:conclusions}, we
summarize our observations and draw conclusions.

\section{Cost functionals subject to linear hyperbolic systems}\label{sec:preliminaries}

\subsection{Problem formulation}
Let $\Omega\subset \sR^d$ ($d=1,2,3$) be an open and bounded domain
with boundary $\Gamma=\partial\Omega$, and let $T>0$. We consider the
linear $n$-dimensional hyperbolic system
\begin{subequations}\label{eq:hyper}
  \begin{alignat}{2}\label{eq:hyper1}
    \qq_t + \nabla\cdot (\FF\qq) &= \ff \qquad &&\text{\ on\ } \Omega\times
    (0,T), \intertext{where, for $(\xx,t)\in\Omega\times (0,T)$, $\qq(\xx,t)\in \sR^n$
      is the vector of state variables and $\ff(\xx,t)\in \sR^n$ is an
      external force. The flux $\FF\qq\in \sR^{n\times d}$ is
      linear in $\qq$ and the divergence operator
      is defined as $\nabla \cdot (\FF\qq) =
      \sum_{i=1}^d {(A_i\qq)}_{x_i}$ with matrix functions $A_i:\Omega\to
      \mathbb{R}^{n\times n}$, where the indices denote partial differentiation with respect to $x_i$.
      Together with~\eqref{eq:hyper1}, we assume the boundary
      and initial conditions}
    B\qq(\xx,t) &= \gg(\xx,t) \qquad &&\xx\in \Gamma,\: t\in
    (0,T),\label{eq:hyper2}\\ \qq(\xx,0) &= \qq_0(\xx) \qquad &&\xx\in
    \Omega.\label{eq:hyper3}
  \end{alignat}
\end{subequations}
Here, $B:\Gamma\to\sR^{l\times n}$ is a matrix function that takes
into account that boundary conditions can only be prescribed on inflow
characteristics.
Under these conditions, \eqref{eq:hyper} has a unique solution $\qq$
in a proper space $Q$~\cite{GustafssonKreissOliger95}.

We target problems, in which the flux, the right hand side, or the
boundary or initial condition data in~\eqref{eq:hyper} depend on
parameters $\cc$ from a space $U$. These parameters can either be
finite-dimensional, i.e., $\cc=(c_1,\ldots,c_k)$ with $k\ge 1$, or
infinite-dimensional, e.g., a function $\cc=\cc(\xx)$. Examples for
functions $\cc$ are material parameters such as the wave speed, or the
right-hand side forcing in~\eqref{eq:hyper1}.

Our main interest are inverse and estimation problems, and optimal control problems
governed by hyperbolic systems of the form~\eqref{eq:hyper}. This leads to
optimization problems of the form
\begin{equation}\label{eq:Ptilde}
  \min_{\cc,\qq} \tilde{\mathcal{J}}(\cc,\qq) \quad \text{subject to }~\eqref{eq:hyper},
\end{equation}
where $\tilde{\mathcal{J}}$ is a cost function that depends on the
parameters $\cc$ and on the state $\qq$.
The parameters $\cc$ may be restricted to an admissible set
$U_{a\!d}\subset U$ for instance to incorporate bound constraints.
If $U_{a\!d}$ is chosen such that for each $\cc\in U_{a\!d}$ the state
equation~\eqref{eq:hyper} admits a unique solution $\qq:=\mathcal
S(\cc)$ (where $\mathcal S$ is the solution operator for the
hyperbolic system), then~\eqref{eq:Ptilde} can be written as an
optimization problem in $\cc$ only, namely
\begin{equation}\label{eq:P}
  \min_{\cc\in U_{a\!d}} \mathcal{J}(\cc) := \tilde{\mathcal{J}}(\cc,\mathcal S(\cc)).
\end{equation}
Existence and (local) uniqueness of solutions to~\eqref{eq:Ptilde} and
\eqref{eq:P} depend on the form of the cost function $\tilde{\mathcal{J}}$,
properties of the solution and parameter spaces and of the hyperbolic
system and have to be studied on a case-to-case basis (we refer, for
instance, to
\cite{BorziSchulz12,Gunzburger03,Lions85,Troltzsch10}). Our main focus
is not the solution of the optimization problem~\eqref{eq:P}, but the
computation of derivatives of $\mathcal J$ with respect to $\cc$, and
the interplay of this derivative computation with the spatial and
temporal discretization of the hyperbolic system
\eqref{eq:hyper}. Gradients (and second derivatives) of
$\mathcal J$ are important to solve~\eqref{eq:P} efficiently,
and can be used for studying parameter sensitivities or quantifying the
uncertainty in the solution of inverse
problems~\cite{Bui-ThanhGhattasMartinEtAl13}.

\subsection{Compatibility of boundary conditions\label{subsec:compatibility}}
To ensure the existence of a solution to the adjoint equation,
compatibility conditions between boundary terms in the cost
function $\tilde{\mathcal J}$, the boundary operator $B$
in~\eqref{eq:hyper2} and the operator $\FF$ in~\eqref{eq:hyper1} must
hold. We consider cost functions of the form
\begin{equation}\label{eq:Jcost}
\tilde{\mathcal{J}}(\cc,\qq) = \int_0^T\!\!\!\int_\Omega j_\Omega(\qq)\,dx\,dt + \int_0^T\!\!\!\int_\Gamma
j_\Gamma(C\qq)\,dx\,dt + \int_\Omega j_T(\qq(T))\,dx,
\end{equation}
where $j_\Omega:\mathbb{R}^n\to\mathbb{R}$,
$j_\Gamma:\mathbb{R}^m\to\mathbb{R}$ and
$J_T:\mathbb{R}^n\to\mathbb{R}$ are differentiable, and
$C:\Gamma\to\sR^{m\times n}$ is a matrix-valued function. We denote
the derivatives of the functional under the integrals by
$j_\Omega'(\cdot)$, $j_\Gamma'(\cdot)$ and $j_T'(\cdot)$. The
Fr\'echet derivative of $\mathcal{J}$ with respect to $\qq$ in a
direction $\tilde \qq$ is given by
\begin{multline}\label{eq:Jcostlin}
\tilde{\mathcal{J}}_{\qq}(\cc,\qq)(\tilde \qq) = \int_0^T\!\!\!\int_\Omega
j'_\Omega(\qq)\tilde \qq\,dx\,dt + \int_0^T\!\!\!\int_\Gamma
j'_\Gamma(C\qq) C\tilde \qq\,dx\,dt\\ + \int_\Omega j'_T(\qq(T))\tilde\qq(T)\,dx.
\end{multline}
The boundary operators $B$ and $C$ must be compatible in the sense
discussed next.  Denoting the outward pointing normal along the boundary $\Gamma$ by
$\nn = {(n_1,\ldots,n_d)}^T$, we use the decomposition
\begin{equation}\label{eq:A}
A := \sum_{i=1}^dn_iA_i = L^T\diag(\lambda_1,\ldots,\lambda_n) L,
\end{equation}
with $L\in \mathbb{R}^{n\times n}$ and $\lambda_1\ge \ldots\ge
\lambda_n$. Note that $L^{-1}=L^T$ if $A$ is symmetric.  The positive
eigenvalues $\lambda_1,\ldots,\lambda_s$, correspond to the $s\ge 0$
incoming characteristics, and the negative eigenvalues
$\lambda_{n-m+1},\ldots,\lambda_n$ to the $m\ge 0$ outgoing
characteristics. Here, we allow for the zero eigenvalues $\lambda_{s+1}=\cdots=\lambda_{n-m}=0$. To ensure well-posedness of the
hyperbolic system $\eqref{eq:hyper}$, the initial values of $\qq$ can
only be specified along incoming characteristics. The first $s$
rows corresponding to incoming characteristics can be identified with
the boundary operator $B$ in~\eqref{eq:hyper2}.  To guarantee
well-posedness of the adjoint equation, $C$ has to be chosen such that
the cost functional $\tilde{\mathcal J}$ only involves boundary
measurements for outgoing characteristics. These correspond to the
rows of $L$ with negative eigenvalues and have to correspond to the
boundary operator $C$. It follows from~\eqref{eq:A} that
\begin{equation}\label{eq:Adecomp}
A =
\begin{bmatrix}
B \\ O \\ C 
\end{bmatrix}^{-1}
\begin{pmatrix}
\lambda_1 \\
& \ddots\\
&& \lambda_n
\end{pmatrix}
\begin{bmatrix}
B \\ O \\ C
\end{bmatrix} = \bar C^T B - \bar B^T C,
\end{equation}
where $O\in \mathbb{R}^{(n-s-m)\times n}$ and $\bar B\in \mathbb{R}^{l\times
  n}$, $\bar C\in \mathbb{R}^{m\times n}$ are derived properly. If
$A$ is symmetric, then $\bar C^T =
B^T\diag(\lambda_1,\ldots,\lambda_S)$, and $\bar B^T =
-C^T\diag(\lambda_{n-m+1},\ldots,\lambda_n)$. As will be shown in
the next section, the matrix $\bar B$ is the boundary condition matrix for the
adjoint equation.
For a discussion of compatibility between the boundary term in a cost
functional and hyperbolic systems in a more general context we refer
to~\cite{GilesPierce97,Hartmann07,AlexeSandu10}.  In the next section,
we formally derive the infinite-dimensional adjoint system and
derivatives of the cost functional $\tilde {\mathcal J}$.

\subsection{Infinite-dimensional derivatives}
For simplicity, we assume that only the
flux $\FF$ (i.e., the matrices $A_1,A_2,A_3$) depend on $\cc$, but
$B$, $C$, $\ff$ and $\qq_0$ do not depend on the parameters $\cc$.
We use the formal Lagrangian method~\cite{Troltzsch10,BorziSchulz12},
in which
we consider $\cc$ and $\qq$ as independent variables and introduce the
Lagrangian function
\begin{align}\label{eq:L}
  \mathscr{L}(\cc,\qq,\pp):= \tilde{\mathcal{J}}(\cc,\qq) +
\int_0^T\!\!\int_\Omega{\left(\qq_t + \nabla\cdot(\FF\qq)-\ff,
\pp\right)}_W\,dx\,dt
\end{align}
with $\pp,\qq\in Q$, where $\qq$ satisfies the boundary and initial
conditions~\eqref{eq:hyper2} and~\eqref{eq:hyper3}, $\pp$ satisfies
homogeneous versions of these conditions, and $\cc\in U_{a\!d}$. Here,
${(\cdot\,,\cdot)}_W$ denotes a $W$-weighted inner product in
$\mathbb{R}^n$, with a symmetric positive definite matrix $W\in
\mathbb{R}^{n\times n}$ (which may depend on $\xx$). The matrix $W$
can be used to make a hyperbolic system symmetric with respect to the
$W$-weighted inner product, as for instance in the acoustic and
coupled elastic-acoustic wave examples discussed in
Sections~\ref{subsec:acoustic} and~\ref{subsec:dgea}. In
particular, this gives the adjoint equation a form very similar to the
state equation.  If boundary conditions depend on the
parameters $\cc$, they must be enforced weakly through a Lagrange
multiplier in the Lagrangian function and cannot be added in the
definition of the solution space for $\qq$. For instance, if the
boundary operator $B=B(\cc)$ depends on $\cc$, the boundary
condition~\eqref{eq:hyper2} must be enforced weakly through a Lagrange
multiplier, amounting to an additional term in the Lagrangian
functional~\eqref{eq:L}.

Following the Lagrangian approach~\cite{Troltzsch10,BorziSchulz12},
the gradient of $\mathcal{J}$
coincides with the gradient of $\mathscr{L}$ with respect to $\cc$,
provided all variations of $\mathscr{L}$ with respect to $\qq$ and
$\pp$ vanish. Requiring that variations  with respect to
$\pp$ vanish, we recover the state
equation. Variations with respect to $\qq$ in directions $\tilde\qq$,
that satisfy homogeneous versions of the initial and boundary
conditions~\eqref{eq:hyper2} and~\eqref{eq:hyper3}, result in
\begin{align*}
\mathscr{L_{\qq}(\cc,\qq,\pp)}(\tilde \qq) %
& = \tilde{\mathcal{J}}_{\qq}(\cc, \qq)(\tilde \qq) -
\int_0^T\!\!\int_\Omega\left(\pp_t,W\tilde\qq\right) +
\left(\FF\tilde\qq,\nabla (W\pp)\right)\,dx\,dt\\
& + \int_\Omega\left(\tilde\qq(T),W\pp(T)\right) +
\int_0^T\!\!\int_{\Gamma}
\left(\nn\cdot\FF\tilde\qq,W\pp\right)\,dx\,dt,
\end{align*}
where we have used integration by parts in time and space. As will be
discussed in Section~\ref{sec:quad}, integration by parts can be
problematic when integrals are approximated using numerical quadrature
and should be avoided to guarantee exact computation of discrete
derivatives.  In this section, we assume continuous functions
$\qq,\pp$ and exact computation of integrals.
Since $\nn\cdot\FF=A$, \eqref{eq:Adecomp} implies that
\begin{equation}\label{eq:compat1}
\int_0^T\!\!\int_{\Gamma}(\nn\cdot\FF\tilde\qq, W\pp) \,dx\,dt=
\int_0^T\!\!\int_{\Gamma}(C\tilde\qq,\bar B W\pp) -
(B\tilde\qq,\bar CW\pp)\,dx\,dt.
\end{equation}
Using the explicit form of the cost given in~\eqref{eq:Jcostlin}, and
that all variations with respect to arbitrary $\tilde \qq$ that
satisfy $B\tilde\qq=0$ must vanish, we obtain
\begin{subequations}\label{eq:adjhyper}
\begin{alignat}{2}
  W\pp_t + \FF^\star\nabla (W\pp) &= j'_\Omega(\qq) \qquad &&\text{\ on\ }
\Omega\times (0,T)\label{eq:adjhyper1} \\
\bar B W\pp(\xx,t) &= -j'_\Gamma(C\qq(\xx,t))
\qquad &&\xx\in \Gamma,\: t\in (0,T),\label{eq:adjhyper2}\\
W\pp(\xx,T) &= -j'_T(\qq(\xx,T)) \qquad &&\xx\in
\Omega.\label{eq:adjhyper3}
\end{alignat}
\end{subequations}
Here, $\FF^\star$ is the adjoint of $\FF$ with respect to the
Euclidean inner product.  Note that the adjoint system
\eqref{eq:adjhyper} is a final value problem and thus is usually solved
backwards in time. Note that, differently from the state system, the
adjoint system is not in conservative form.

Next we compute variations of $\mathscr L$ with respect to the
parameters $\cc$ and obtain for variations $\tilde \cc$ that
\begin{subequations}\label{eq:Jgrad}
\begin{align}
\mathscr L_\cc(\cc,\qq,\pp)(\tilde \cc) &= \int_0^T\!\!\!\int_\Omega
{\left(\nabla\cdot(\FF_\cc(\tilde \cc)\qq), \pp\right)}_W =
 \int_0^T\!\!\!\int_\Omega
\sum_{i=1}^d{\left({({A_i}_\cc(\tilde\cc) \qq)}_{x_i}, \pp\right)}_W \label{eq:Jgrad1}%
\end{align}
\end{subequations}
Since $\qq$ and $\pp$ are assumed to solve the
state and adjoint system,
\begin{equation}
\mathcal J(\cc)(\tilde\cc) = \mathscr L_\cc(\cc,\qq,\pp)(\tilde \cc),
\ \text{where}\ \qq\ \text{solves~\eqref{eq:hyper}}\ \text{and}\ \pp
\ \text{solves~\eqref{eq:adjhyper}}.
\end{equation}
Next, we present the dG discretization of the hyperbolic system
\eqref{eq:hyper} and discuss the interaction between discretization
and the computation of derivatives.

\subsection{Discontinuous Galerkin discretization}\label{subsec:dG}
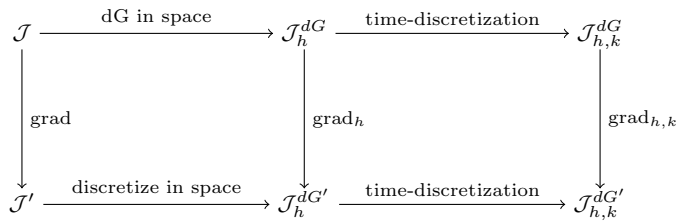
\begin{figure}\centering
\begin{tikzpicture}[description/.style={fill=white,inner sep=2pt}]
\matrix (m) [matrix of math nodes, row sep=5em,
column sep=10em, text height=1.5ex, text depth=0.25ex]
{\mathcal{J} &
\mathcal{J}^{dG}_h & \mathcal{J}^{dG}_{h,k}  \\[3ex]
\mathcal{J}'  & \mathcal{J}^{dG'}_{h} &  \mathcal{J}^{dG'}_{h,k} \\};
\path[->,font=\scriptsize]
(m-1-1) edge node[auto] {\text{grad}} (m-2-1)
(m-1-2) edge node[auto] {\text{grad$_h$}} (m-2-2)
(m-1-3) edge node[auto] {\text{grad$_{h,k}$}} (m-2-3)
(m-1-1) edge node[auto] {\text{dG in space}} (m-1-2)
(m-1-2) edge node[auto] {\text{time-discretization}} (m-1-3)
(m-2-1) edge node[auto] {\text{discretize in space}} (m-2-2)
(m-2-2) edge node[auto] {\text{time-discretization}} (m-2-3);
\end{tikzpicture}
\caption{Sketch to illustrate the relation between discretization of
  the problem (horizontal arrows, upper row), computation of the
  gradient with respect to the parameters $\cc$ (vertical arrows) and
  discretization of the gradient (horizontal arrows, lower row). The
  problem discretization (upper row) requires discretization of the
  state equation and of the cost functional in space (upper left
  horizontal arrow) and in time (upper right horizontal arrow). The
  vertical arrows represent the Fr\`echet derivatives of $\mathcal J$
  (left), of the semidiscrete cost $\mathcal{J}^{dG}_h$ (middle) and
  the fully discrete cost $\mathcal{J}^{dG}_{h,k}$ (right).  The
  discretization of the gradient (bottom row) requires space (left
  arrow) and time (right arrow) discretization of the state equation,
  the adjoint equation and the expression for the gradient. Most of
  our derivations follow the fully discrete approach, i.e., the upper
  row and right arrows; The resulting discrete expressions are then
  interpreted as discretizations of the corresponding continuous
  equations derived by following the vertical left
  arrow.\label{fig:diagram}}
\end{figure}

For the spatial discretization of hyperbolic systems such as
\eqref{eq:hyper}, the discontinuous Galerkin (dG) method has proven to
be a favorable choice.
In the dG method, we divide the domain $\Omega$ into disjoint elements
$\De$, and use polynomials to approximate $\qq$ on each element
$\De$. The resulting approximation space is denoted by $Q^h$, and
elements $\qq_h\in Q^h$ are polynomial on each element, and
discontinuous across elements. Using test functions $\pp_h\in Q^h$,
dG discretization in space implies that for each element $\De$
\begin{equation}\label{eq:dGweak}%
  \begin{split}
    \int_{\De}(\frac \partial {\partial t}\qq_h, W\pp_h) & -
    (\FF\qq_h, \nabla (W\pp_h)) \,d\xx +\\
    &\int_{\Gamma^e}\nn^-\cdot({(\FF\qq_h)}^\dagger, W\pp_h^-)\,d\xx =
    \int_{\De}(\ff,W\pp_h)\,d\xx
  \end{split}
\end{equation}
for all times $t\in (0,T)$. Here, $(\cdot\,,\cdot)$ is the inner
product in $\sR^n$ and $\sR^{d\times n}$ and the symmetric and
positive definite matrix $W$ acts as a weighting matrix in this inner
product. Furthermore, ${(\FF\qq_h)}^\dagger$ is the numerical flux, which
connects adjacent elements. The superscript ``$-$'' denotes that the
inward values are chosen on $\Gamma^e$, i.e., the values of the
approximation on $\Omega^e$; the superscript ``$+$'' denotes that the
outwards values are chosen, i.e., the values of an element
$\Omega^{e'}$ that is adjacent to $\Omega^e$ along the shared boundary $\Gamma^e$.
Here, $\nn^-$ is the outward pointing normal on element
$\Omega^e$.
The
formulation~\eqref{eq:dGweak} is often referred to as the \emph{weak form} of the dG
discretization~\cite{HesthavenWarburton08,Kopriva09}. The
corresponding \emph{strong form} dG discretization is obtained by
element-wise integration by parts in space in~\eqref{eq:dGweak},
resulting in
\begin{equation}\label{eq:dGstrong}%
  \begin{split}
    \int_{\De}(\frac\partial{\partial t}\qq_h, W\pp_h) & +
    (\nabla \cdot\FF\qq_h, W\pp_h) \,d\xx -\\
    &\int_{\Gamma^e}\nn^-\cdot(\FF\qq_h^- - {(\FF\qq_h)}^\dagger,
    W\pp_h^-)\,d\xx = \int_{\De}(\ff,W\pp_h)\,d\xx
  \end{split}
\end{equation}
for all $t\in (0,T)$.
To find a solution to the optimization problem~\eqref{eq:P},
derivatives of $\mathcal J$ with respect to the parameters $\cc$ must
be computed. There are two choices for computing derivatives, namely
deriving expressions for the derivatives of the continuous
problem~\eqref{eq:P}, and then discretizing these equations, or first
discretizing the problem and then computing derivatives of this fully
discrete problem.  If the latter approach is taken, i.e., the discrete
adjoints are computed, the question arises weather the discrete
adjoint equation is an approximation of the continuous adjoint and if
the discrete adjoint equation is again a dG discretization. Moreover,
what are the consequences of choosing the weak or the strong
form~\eqref{eq:dGweak} or~\eqref{eq:dGstrong}?
A sketch for the different combinations of discretization
and computation of derivatives is also shown in
Figure~\ref{fig:diagram}.

\subsection{Influence of numerical quadrature}
\label{sec:quad}
In a numerical implementation, integrals are often approximated using
numerical quadrature. A fully discrete approach has to take into
account the resulting quadrature error; in particular, integration by
parts can incur an error in combination with numerical
quadrature. Below, we first discuss implications of numerical
quadrature in space and then comment on numerical integration
in time. In our example problems in Section~\ref{sec:examples}, we use
integral symbols to denote integration in space and time, but do not
assume exact integration. In particular, we avoid integration by parts
or highlight when integration by parts is used.

\subsubsection{Numerical integration in space}
\label{subsec:quadrature-in-space}
The weak form~\eqref{eq:dGweak} and the strong form~\eqref{eq:dGstrong}
of dG are equivalent provided integrals are computed exactly and, as a
consequence, integration by parts does not result in numerical
error. If numerical quadrature is used, these forms are only
numerically equivalent under certain
conditions~\cite{KoprivaGassner10}; in general, they are different.
To compute fully discrete gradients, we thus avoid integration by
parts in space whenever possible. As a consequence, if the weak form
for the state equation is used, the adjoint equation is in strong
form; this is illustrated and further discussed in
Section~\ref{sec:examples}.

\subsubsection{Numerical integration in time}\label{subsec:quadrature-in-time}

We use a method-of-lines approach,
that is, from the dG discretization in space we obtain a
continuous-in-time system of ordinary differential equations (ODEs),
which is then discretized by a Runge--Kutta method.  Discrete adjoint
equations and the convergence to their continuous counterparts for
systems of ODEs discretized by Runge--Kutta methods have been studied,
for instance in~\cite{Hager00, Walther07}. For the time-discretization
we build on these results.

An alternative approach to discretize in time is using a finite
element method for the time discretization, which allows a fully
variational formulation of the problem in space-time; we refer, for
instance to~\cite{BeckerMeidnerVexler07} for this approach applied to
parabolic optimization problems.

In both approaches, the computation of derivatives
requires the entire time history of
both the state and the adjoint solutions. For realistic application
problems, storing this entire time history is infeasible, and storage
reduction techniques, also known as checkpointing strategies have to
be employed.  These methods allow to trade storage against computation
time by storing the state solution only at certain time instances, and
then recomputing it as needed when solving the adjoint
equation and computing the
gradient~\cite{BeckerMeidnerVexler07,GriewankWalther00,GriewankWalther08}.

\section{Example problems}\label{sec:examples}
The purpose of this section is to illustrate the issues discussed in
the previous sections on example problems.  We present examples with
increasing complexity and with parameters entering differently in the
hyperbolic systems. First, in Section~\ref{subsec:advection}, we
derive expressions for derivatives of a functional that depends on the
solution of the one-dimensional advection equation. Since linear
conservation laws can be transformed in systems of advection
equations, we provide extensive details for this example. In
particular, we discuss different numerical fluxes. In
Section~\ref{subsec:acoustic}, we compute expressions for the
derivatives of a functional with respect to the local wave speed in an
acoustic wave equation. This is followed by examples in which we
compute derivatives with respect to a boundary forcing in Maxwell's
equation (Section~\ref{subsec:maxwell}) and derivatives with respect
to the primary and secondary wave speeds in the coupled
acoustic-elastic wave equation (Section~\ref{subsec:dgea}).  At the
end of each example, we summarize our observations in remarks.

Throughout this section, we use the dG discretization introduced in
Section~\ref{subsec:dG} and denote the finite dimensional dG solution
spaces by $P^h$ and $Q^h$. These spaces do not include the boundary
conditions which are usually imposed weakly through the numerical flux
in the dG method, and they also do not include initial/final time
conditions, which we specify explicitly. Functions in these dG spaces
are smooth (for instance polynomials) on each element $\Omega^e$ and
discontinuous across the element boundaries $\partial\De$.  As before,
for each element we denote the inward value with a superscript ``-''
and the outward value with superscript ``+''. We use the
index $h$ to denote discretized fields and denote by $\jump{\cdot}$
the jump, by $\diff{\cdot}$ the difference and by $\mean{\cdot}$ the
mean value at an element interface $\partial\De$. To be precise, for a
scalar dG-function $u_h$, these are defined as $\jump{u_h}=u_h^-\nn^-
+ u_h^+\nn^+$, $\diff{u_h}=u_h^--u_h^+$ and
$\mean{u_h}=(u_h^-+u_h^+)/2$.  Likewise, for a vector $\vv_h$ we
define $\jump{\vv_h}=\vv_h^-\cdot\nn^- + \vv_h^+\cdot\nn^+$,
$\diff{\vv_h}=\vv_h^--\vv_h^+$ and $\mean{\vv_h}=(\vv_h^-+\vv_h^+)/2$
and for a second-order tensor $\SS_h$ we have
$\jump{\SS_h}=\SS_h^-\nn^-+\SS_h^+\nn^+$.  The domain boundary $\partial\Omega$
is denoted by $\Gamma$.

Throughout this section, we use the usual symbol to denote integrals,
but we do not assume exact quadrature. Rather, integration can be
replaced by a numerical quadrature rule and, as a consequence, the
integration by parts formula does not hold exactly.  To avoid
numerical errors when using numerical quadrature, we thus avoid
integration by parts in space or point out when integration by parts
is used. Since our focus is on the spatial dG discretization, we do
not discretize the problem in time and assume exact integration in
time.

\subsection{One-dimensional advection equation}\label{subsec:advection}
We consider the one-dimensional advection equation on the spatial
domain $\Omega = (x_l,x_r)\subset \mathbb{R}$. We assume a spatially
varying, continuous positive advection velocity
$a(x)\ge a_0>0$ for $x\in\Omega$ and a forcing $f(x,t)$ for $(x,t)\in
\Omega\times (0,T)$. The advection equation written in conservative
form is given by
\begin{subequations}
  \begin{alignat}{2}\label{eq:adv}
    u_t + {(au)}_x &= f \quad &&\text{\ on\ } \Omega\times (0,T),
    \intertext{with the initial condition}
    u(x,0) &= u_0(x) \quad &&\text{\ for\ } x\in \Omega,\\
    \intertext{and, since $a>0$, the inflow boundary is $\Gamma_l:=\{x_l\}$, where we assume}
    u(x_l,t) &= u_l(t) \quad &&\text{\ for\ } t\in (0,T).
  \end{alignat}
\end{subequations}
The discontinuous Galerkin (dG) method for the numerical solution of
\eqref{eq:adv} in strong form is: Find $u_h\in P^h$ with
$u_h(x,0)=u_0(x), x\in \Omega$ such that for all test functions
$p_h\in P^h$ holds
\begin{equation}\label{adv:strong}
\int_\Omega (u_{h,t} + {(au_h)}_x - f)p_h \,dx = \sum_e\int_{\partial\De} n^- \left( au_h^- -
{(au_h)}^\dagger \right)p_h^- \,dx
\end{equation}
for all $t\in (0,T)$. Here, $a\in U$ can be an infinite-dimensional continuous
function, or a finite element function.
For $\alpha\in [0,1]$, the numerical flux on the
boundary is replaced by the numerical flux, ${(au_h)}^\dagger$, given
by
\begin{equation}\label{adv:flux}
{(au_h)}^\dagger  =  a \mean{u_h} + \frac 12 \abs{a}(1-\alpha) \jump{u_h}.
\end{equation}
This is a central flux for $\alpha=1$, and an upwind flux for
$\alpha=0$.  Note that in one spatial dimension, the outward normal
$n$ is $-1$ and $+1$ on the left and right side of $\De$,
respectively. Since $a$ is assumed to be continuous on $\Omega$, we
have $a^- = a^+:=a$. Moreover, since $a$ is positive, we neglect the
absolute value in the following.

As is standard practice~\cite{HesthavenWarburton08,Kopriva09} we
incorporate the boundary conditions weakly through the numerical flux
by choosing the ``outside'' values $u_h^+$ as
$ u_h^+ = u_l$ for $x=x_l$ and $u_h^+ = u_h^-$ for $x=x_r$
for the computation of ${(au_h)}^\dagger $. This implies that
$\jump{u_h}=n^-(u_h^--u_l)$ on $\Gamma_l$ and $\jump{u_h}=0$ on
the outflow boundary $\Gamma_r:=\{x_r\}$.
For completeness, we also provide the dG discretization of
\eqref{eq:adv} in weak form: Find $u_h\in P^h$ with $u_h(x,0)=u_0(x),
x\in \Omega$ such that for all $p_h\in P^h$ holds
\begin{equation*}
\int_\Omega (u_{h,t} - f)p_h - au_h p_{h,x} \,dx = -\sum_e\int_{\partial\De} n^-
{(au_h)}^\dagger p_h^- \,dx %
\end{equation*}
for all $t\in (0,T)$, which is found
by element-wise integration by parts in~\eqref{adv:strong}.

Adding the boundary contributions in~\eqref{adv:strong} from adjacent
elements $\De$ and $\Dep$ to their shared point $\partial\De\cap
\partial\Dep$, we obtain
\begin{equation*}
\frac 12 \abs{a} (\alpha-1)\jump{u_h}\jump{p_h} + a\jump{u_h}\mean{p_h}.
\end{equation*}
Thus, \eqref{adv:strong} can also be written as
\begin{multline}\label{adv:strong1}
  \int_\Omega (u_{h,t} + {(au_h)}_x - f)p_h \,dx = \sum_e\int_{\partial\De\setminus\Gamma} n^-
  \left(a\mean{p_h} + {a}\frac 12 (\alpha-1)\jump{p_h}\right)u_h^- \,dx\\
  + \int_{\Gamma_l} n^- \left(a\frac12 p_h^- + {a}\frac 12
    (\alpha-1)n^-p_h^-\right)(u_h^--u_l) \,dx.
\end{multline}

We consider an objective functional for the advection velocity $a\in U$ given by
\begin{equation}\label{eq:cost}
\mathcal J(a):=\tilde{\mathcal{J}}(a,u_h) = \int_0^T\!\!\int_\Omega
j_\Omega(u_h)\,dx\,dt + \int_0^T\!\!\int_{\Gamma} j_\Gamma(u_h)\,dx\,dt
+ \int_\Omega r(a)\,dx,
\end{equation}
with differentiable functions $j_\Omega:\Omega\to\mathbb{R}$,
$j_\Gamma:\Gamma\to\mathbb{R}$ and $r:\Omega\to\mathbb{R}$.  For
illustration purposes, we define the boundary term in the cost on both,
the inflow and the outflow part of the boundary, and comment on the
consequences in Remark~\ref{rem:adv1}.
To derive the discrete gradient of $\mathcal{J}$, we use the
Lagrangian function $\mathcal L: U\times P^h\times P^h\rightarrow
\mathbb R$, which combines the cost~\eqref{eq:cost} with the dG
discretization~\eqref{adv:strong1}:
\begin{align*}
\mathcal{L}(a,u_h,p_h) :=& \tilde{\mathcal{J}}(a,u_h) + \int_0^T\!\!\int_\Omega (u_{h,t} + {(au_h)}_x - f)p_h \,dx\,dt \\&-
\sum_e\int_0^T\!\!\int_{\partial\De\setminus \Gamma} n^-
\left(a\mean{p_h} + {a}\frac 12 (\alpha-1)\jump{p_h}\right)u_h^- \,dx\,dt\\
&- \int_0^T\int_{\Gamma_l} n^- \left(a\frac12 p_h^- + {a}\frac 12
  (\alpha-1)n^-p_h^- \right)(u_h^--u_l) \,dx\,dt.
\end{align*}
By requiring that all variations with respect to $p_h$ vanish, we recover the state
equation~\eqref{adv:strong1}. Variations with respect to $u_h$ in a
direction $\tilde u_h$, which satisfies the homogeneous initial
conditions $\tilde u_h(x,0) = 0$ for $x\in\Omega$ result in
\begin{align*}
&\!\!\!\!\mathcal L_{u_h}(a,u_h,p_h)(\tilde u_h)\\
 =& \int_0^T\!\!\int_\Omega j'_\Omega(u_h)\tilde u_h
-p_{h,t} \tilde u_h + {(a\tilde u_h)}_x p_h \,dx\,dt + \int_\Omega p_h(x,T)\tilde u_h(x,T)\,dx
 \\&- \sum_e\int_0^T\!\!\int_{\partial\De\setminus\Gamma}\!\!\!
n^- \left(a\mean{p_h} + {a}\frac 12 (\alpha-1)\jump{p_h}\right)\tilde u_h^-
\,dx\,dt + \int_0^T\!\!\int_{\Gamma} j'_\Gamma(u_h)\tilde u_h\,dx\,dt\\
&- \int_0^T\!\!\int_{\Gamma_l}
n^- \left(a\frac12 p_h^- + {a}\frac 12 (\alpha-1)n^-p_h^-\right)\tilde u_h^-\,dx\,dt
\end{align*}
Since we require that arbitrary variations with respect to $u_h$ must
vanish, $p_h$ has to satisfy $p_h(x,T) = 0$ for all $x\in\Omega$, and
\begin{subequations}\label{adv:adjweak}
  \begin{multline}\label{adv:adjweak1}
      \int_\Omega -p_{h,t} \tilde u_h + {(a\tilde u_h)}_x p_h + j_\Omega'(u_h)\tilde u_h \,dx =
      \sum_e\int_{\partial\De\setminus\Gamma_l} n^- {(ap_h)}^\dagger \tilde u_h^- \,dx \\+
      \int_{\Gamma_l} n^- \left(\frac{a(2-\alpha)}{2} p_h^- + j_\Gamma'(u_h)\right)\tilde u_h^-\,dx
  \end{multline}
 for all $t\in (0,T)$ and for all $\tilde u_h$, with the adjoint flux
  \begin{equation}\label{adv:adjflux}
    {(ap_h)}^\dagger:=  a\mean{p_h} + \frac{1}{2} {a}(\alpha-1)\jump{p_h}
  \end{equation}
  and with
  \begin{equation}\label{adv:adjfluxbc}
    p_h^+ :=  -\frac{j_\Gamma'(u_h)}{a(1-\frac\alpha2)} -
    \frac{\alpha}{2-\alpha}p_h^- \quad\text{on}\ \Gamma_r.
  \end{equation}
\end{subequations}
The outside value $p_h^+$ in~\eqref{adv:adjfluxbc} is computed such that
$n^-{(ap_h)}^\dagger = j'_\Gamma(u_h)$ on $\Gamma_r$.  The equations
\eqref{adv:adjweak} are the weak form of a discontinuous Galerkin
discretization of the adjoint equation, with flux given by
\eqref{adv:adjflux}.
An element-wise integration by parts in space in~\eqref{adv:adjweak1},
results in the corresponding strong form of the discrete adjoint
equation:
\begin{equation}\label{adv:adjstrong}
  \begin{split}
  \int_\Omega \big(-p_{h,t} - a p_{h,x} + & j_\Omega'(u_h)\big)\tilde u_h \,dx =
  -\sum_e\int_{\partial\De\setminus\Gamma_l} n^-\left(ap_h^- - {(ap_h)}^\dagger\right)
  \tilde u_h^- \,dx\\
  & + \int_{\Gamma_l} n^- \left(-\frac{a\alpha}{2} p_h^- + j_\Gamma'(u_h)\right)\tilde u_h^-\,dx
  \end{split}
\end{equation}
Note that this integration by parts is not exact if numerical quadrature
is used. It can be avoided if the dG weak form of the adjoint equation
\eqref{adv:adjweak} is implemented directly.

Provided $u_h$ and $p_h$ are solutions to the state and adjoint
equations, respectively, the gradient of $\mathcal J$ with respect to
$a$ is found by taking variations of the Lagrangian with respect to
$a$ in a direction $\tilde a$:
\begin{align*}
\mathcal L_a(a,u_h,p_h)(\tilde a) =& \int_\Omega r'(a)\tilde a\,dx +
\int_0^T\!\!\int_\Omega {(\tilde a u_h)}_x p_h \,dx\,dt \\&-
\sum_e\int_0^T\!\!\int_{\partial\De\setminus\Gamma} n^- \tilde a
\left(\frac 12 (\alpha-1)\jump{p_h} + \mean{p_h} \right) u_h^- \,dx\,dt\\ &-
\int_0^T\!\!\int_{\Gamma_l} n^- \tilde a\left(\frac 12
(\alpha-1)n^-p_h^- + \frac12 p_h^- \right)(u_h^--u_l) \,dx\,dt.
\end{align*}
Thus, the gradient of $\mathcal J$ is given by
\begin{equation}\label{adv:grad}
{\mathcal J}'(a)(\tilde a) = \int_\Omega r'(a)\tilde a \,dx +
\int_0^T\!\!\int_\Omega G\tilde a \,dx\,dt +
\sum_{\mathsf f} \int_0^T\!\!\int_{\mathsf f} g\tilde a\,dx\,dt,
\end{equation}
where $G$ is defined for each element $\De$, and $g$ for each
inter-element face as follows:
\begin{align*}%
G\tilde a = p_h{(\tilde a u_h)}_x \quad \text{and}\quad g =
\begin{cases}
\frac 12 (\alpha-1)\jump{u_h}\jump{p_h} + \jump{u_h}\mean{p_h} &\ \text{if}\ \mathsf f\not\subset \Gamma,\\
\frac 12 \alpha p_h^-(u_h^--u_l) &\ \text{if}\ \mathsf f\subset \Gamma_l,\\
0  &\ \text{if}\ \mathsf f\subset \Gamma_r.
\end{cases}
\end{align*}
Note that the element boundary jump terms arising in $\mathcal J'(a)$
are a consequence of using the dG method to discretize the state
equation. While in general these terms do not vanish, they become
small as the discretization resolves the continuous state and adjoint
variables. However, these terms must be taken into account to compute
discretely exact gradients.  We continue with a series of remarks:

\begin{remark}\label{rem:adv1}
The boundary conditions for the adjoint variable $p_h$ that are weakly
imposed through the adjoint numerical flux are the Dirichlet condition
$a p_h(x_r) = -j'_\Gamma(u)$ at $\Gamma_r$, which, due to the sign
change for the advection term in~\eqref{adv:adjweak} and
\eqref{adv:adjstrong}, is an inflow boundary for the adjoint equation.
At the (adjoint) outflow boundary $\Gamma_l$,
the adjoint scheme can only be stable if
${j_\Gamma}|_{\Gamma_l}\equiv 0$. This corresponds to the discussion
from Section~\ref{subsec:compatibility} on the compatibility of boundary operators.
The discrete adjoint scheme is consistent (in the sense that the
continuous adjoint variable $p$ satisfies the discrete adjoint
equation~\eqref{adv:adjweak}) when $\alpha=0$, i.e., for a dG
scheme based on an upwind numerical flux. Thus, only dG
discretizations based on upwind fluxes at the boundary can be used in adjoint
calculus. Hence, in the following we restrict ourselves to upwind
fluxes.
\end{remark}
\begin{remark}\label{rem:adv2}
While the dG discretization of the state equation is in conservative
form, the discrete adjoint equation is not.  Moreover, using dG method
in strong form for the state system, the adjoint system is naturally dG method in weak form (see
\eqref{adv:adjweak}), and element-wise integration by parts is
necessary to find the adjoint in strong form
\eqref{adv:adjstrong}. Vice versa, using the weak form of dG for the
state equation, the adjoint equation is naturally in strong
form. These two forms can be numerically different if the integrals
are approximated through a quadrature rule for which
integration by parts does not hold exactly. In this case, integration by parts should be avoided to obtain exact
discrete derivatives.
\end{remark}
\begin{remark}\label{rem:adv3}
The numerical fluxes~\eqref{adv:flux} and~\eqref{adv:adjflux} differ
by the sign in the upwinding term only. Thus, an upwinding flux for
the state equation becomes a downwinding flux for the adjoint
equation. This is natural since the advection velocity for the adjoint
equation is $-a$, which makes the adjoint numerical flux an upwind
flux for the adjoint equation.
\end{remark}

\subsection{Acoustic wave equation}\label{subsec:acoustic}
Next, we derive expressions for the discrete gradient with respect to
the local wave speed in the acoustic wave equation. This is important,
for instance, in seismic inversion using full wave forms
\cite{Fichtner11, LekicRomanowicz11,
  EpanomeritakisAkccelikGhattasEtAl08, PeterKomatitschLuoEtAl11}.
If the dG method is used to discretize the wave equation (as, e.g., in
\cite{Bui-ThanhBursteddeGhattasEtAl12, Bui-ThanhGhattasMartinEtAl13,
  CollisObervanBloemenWaanders10}), the question on the proper
discretization of the adjoint equation and of the expressions for the
derivatives arises. Note that in Section~\ref{subsec:dgea} we
present the discrete derivatives with respect to the
(possibly discontinuous) primary and secondary wave speeds in the
coupled acoustic-elastic wave equation, generalizing the results
presented in this section. However, for better readability we choose
to present this simpler case first and then present the results for
the coupled acoustic-elastic equation in compact form in
Section~\ref{subsec:dgea}.

We consider the acoustic wave equation written as first-order system
as follows:
\begin{subequations}\label{eq:acoustic}
  \begin{alignat}{2}
    e_t - \nabla \cdot \vv &= 0 \quad&\text{\ on\ } \Omega\times (0,T),\label{eq:acoustic1}\\
    \rho \vv_t - \nabla (\lambda e) &= \ff \quad&\text{\ on\ } \Omega\times (0,T),\label{eq:acoustic2}
  \end{alignat}
  where $\vv$ is the velocity, $e$ the dilatation (trace of the strain
  tensor), $\rho=\rho(\point{x})$ is the mass density, and $\lambda =
  c^2\rho$, where $c(\point{x})$ denotes the wave speed.  Together
  with~\eqref{eq:acoustic1} and~\eqref{eq:acoustic2}, we assume the
  initial conditions
  \begin{align}\label{eq:acoustic3}
    e(\xx,0) = e_0(\xx), \:\vv(\xx,0) = \vv_0(\xx) \quad \text{\ for\ } \xx\in \Omega,
  \end{align}
  and the boundary conditions
  \begin{align}
    e(\xx,t) = e_{\text{bc}}(\xx,t),\:\: \vv(\xx,t) = \vv_\text{bc}(\xx,t)\quad \text{\ for\ } (\xx,t)\in \Gamma\times(0,T).
  \end{align}
\end{subequations}
Note that the dG method discussed below
uses an upwind numerical flux, and thus the boundary conditions
\eqref{eq:acoustic} are automatically only imposed at inflow
boundaries. Through proper choice of $e_{\text{bc}}$ and
$\vv_\text{bc}$ classical wave equation boundary conditions can be
imposed, e.g., \cite{WilcoxStadlerBursteddeEtAl10, FengTengChen07}.

 The choice of the dilatation $e$ together with the
velocity $\vv$ in the first order system formulation is motivated from
the strain-velocity formulation used for the coupled elastic and
acoustic wave equation in Section~\ref{subsec:dgea}.
To write~\eqref{eq:acoustic1} and~\eqref{eq:acoustic2} in second-order form, we
define the pressure as
$p=-\lambda e$ and obtain the 
pressure-velocity form as
\begin{subequations}\label{eq:acoustic_p}
  \begin{alignat*}{2}
    p_t + \lambda \nabla \cdot \vv &= 0 \quad&\text{\ on\ } \Omega\times (0,T),\\%\label{eq:acoustic_p1}\\
    \rho \vv_t + \nabla p &= \ff \quad&\text{\ on\ } \Omega\times (0,T),%
  \end{alignat*}
\end{subequations}
which is equivalent to the second-order formulation
\begin{equation}
p_{tt} = \lambda\nabla\cdot\left(\frac 1\rho\nabla p \right) - \lambda\nabla\cdot\left(\frac{1}{\rho}\ff\right)\quad \text{\ on\ } \Omega\times (0,T).
\end{equation}
The strong form dG discretization of~\eqref{eq:acoustic} is:
Find $(e_h,\vv_h)\in P^h\times Q^h$ satisfying the initial conditions
\eqref{eq:acoustic3} such that for all test functions $(h_h,\ww_h)\in
P^h\times Q^h$ and for all $t\in (0,T)$ holds:
\begin{equation}\label{acoustic:strong}
  \begin{split}
    &\int_\Omega (e_{t,h} - \nabla\cdot \vv_h) \lambda h_h \,d\xx
    + \int_\Omega (\rho\vv_{t,h} - \nabla (\lambda e_h) -\ff)\cdot\ww_h\,d\xx\\
    &  = -\sum_e\int_{\partial\De} \nn^-\cdot\left( \vv_h^- - \vv_h^\dagger\right)\lambda
    h_h^- + \left( {(\lambda e_h)}^- - {(\lambda
    e_h)}^\dagger\right) \nn^-\cdot\ww_h^- \,d\xx.
  \end{split}
\end{equation}
Note that above, the inner product used for~\eqref{eq:acoustic1} is
weighted by $\lambda$, which makes the first-order form of the wave
equation a symmetric operator and also allows for a natural
interpretation of the adjoint variables, as shown below. Assuming that
$c$ and $\rho$ are continuous, we obtain the upwind numerical fluxes:
\begin{subequations}\label{acoustic:flux}
  \begin{align}
    \nn^-\cdot\vv_h^\dagger  &= \nn^-\cdot\mean{\vv_h} - \frac c2\diff{e_h},   \\
    {(\lambda e_h)}^\dagger &= \lambda \mean{e_h} - \frac{\rho
      c}{2}\jump{\vv_h}.
  \end{align}
\end{subequations}
Adding the boundary contributions from two adjacent elements $\De$ and
$\Dep$ in~\eqref{acoustic:strong} to a shared edge (in 2D) or face (in
3D), one obtains
\begin{equation}\label{adv:boundary}
\lambda \jump{\vv_h}\mean{h_h} + \frac{c\lambda}{2}\diff{e_h}\diff{h_h} + \frac{\rho
c}{2} \jump{\vv_h}\jump{\ww_h} + \lambda\jump{e_h}\cdot\mean{\ww_h}.
\end{equation}
We compute the discrete gradient with respect to the wave speed
$c$ for a cost functional of the form
\begin{equation}\label{eq:acoustic-cost}
\tilde{\mathcal{J}}(c,\vv_h) = \int_0^T\!\!\int_\Omega
j_\Omega(\vv_h)\,d\xx\,dt + \int_0^T\!\!\int_\Gamma j_\Gamma(\vv_h)\,d\xx\,dt
+ \int_\Omega r(c)\,d\xx,
\end{equation}
with differentiable functions $j_\Omega:\Omega\to\mathbb{R}$,
$j_\Gamma:\Gamma\to\mathbb{R}$ and $r:\Omega\to\mathbb{R}$. To ensure
compatibility as discussed in Section~\ref{subsec:compatibility}, the
boundary term $j_\Gamma$ in~\eqref{eq:acoustic-cost} can only involve
outgoing characteristics. We
introduce the Lagrangian function, use that all its variations with
respect to $\vv$ and $e$ must vanish, and integrate by parts in time
$t$, resulting in the following adjoint equation: Find $(h_h,\ww_h)\in
P^h\times Q^h$ satisfying the finial time conditions $h(\xx,T)=0$,
$\ww(\xx,T)=0$ for $\xx\in \Omega$, such that for all
  test functions $(\tilde e_h,\tilde \vv_h)\in P^h\times Q^h$ and for
  all $t\in (0,T)$ holds:
\begin{equation}\label{acoustic:adjweak}
  \begin{split}
    &\int_\Omega -\tilde e_h\lambda h_{h,t} - \nabla\cdot \tilde \vv_h \lambda
    h_h -\rho\tilde\vv_h\cdot \ww_{h,t} - \nabla (\lambda \tilde e_h)\cdot\ww_h +
    j_\Omega'(\vv_h)\cdot\tilde \vv_h\,d\xx \\
    &  = -\sum_e\int_{\partial\De} \left(\nn^-\cdot \ww_h^\dagger\right) \lambda
    \tilde e_h^- + {(\lambda h_h)}^\dagger \nn^-\cdot\tilde \vv_h^- \,d\xx -
    \int_\Gamma j_\Gamma'(\vv_h)\cdot\tilde \vv_h\,d\xx,
  \end{split}
\end{equation}
where the adjoint numerical fluxes are given by
\begin{subequations}\label{acoustic:adjflux}
\begin{align}
  \nn^-\cdot\ww_h^\dagger &= \nn^-\cdot\mean{\ww_h} + \frac c2 \diff{h_h},   \\
  {(\lambda h_h)}^\dagger   &= \lambda\mean{h_h} + \frac{\rho c}{2} \jump{\ww_h}.
\end{align}
\end{subequations}
Note that~\eqref{acoustic:adjweak} is the weak form of the dG
discretization for an acoustic wave equation, solved backwards in
time. This is a consequence of the symmetry of the differential
operator in the acoustic wave equation, when considered in the
appropriate inner product.  Comparing~\eqref{acoustic:adjflux}
and~\eqref{acoustic:flux} shows that the
adjoint numerical flux~\eqref{acoustic:adjflux} is the downwind flux
in the adjoint variables for the adjoint wave equation. The strong dG
form corresponding to~\eqref{acoustic:adjweak} can be obtained by
element-wise integration in parts in space.

Finally, we present expressions for the derivative of $\mathcal J$
with respect to the wave speed $c$, which are found as variations
of the Lagrangian with respect to $c$. This results in
\begin{equation}\label{acoustic:grad_gen}
{\mathcal J}'(c)(\tilde c) = \int_\Omega r'(c)\tilde c \,d\xx +
\int_0^T\!\!\int_\Omega G\tilde c \,d\xx\,dt +
\sum_{\mathsf f} \int_0^T\!\!\int_{\mathsf f} g\tilde c\,d\xx\,dt,
\end{equation}
where $G$ is defined on each element $\De$, and $g$ for each
inter-element face $\mathsf f$ as follows:
\begin{subequations}\label{acoustic:grad}
  \begin{align}\label{acoustic:grad1}
    G \tilde c &= -2\nabla(\rho c \tilde c e_h)\cdot \ww_h
     +2\rho c\tilde c (e_{h,t} - \nabla\cdot\vv_h)h_h, \\
    g &= 2\rho c \jump{\vv_h}\mean{h_h} + \frac 32 \rho c^2 \diff{e_h}\diff{h_h} +
    \frac 12 \rho \jump{\vv_h}\jump{\ww_h} + 2 \rho
    c\jump{e_h}\cdot\mean{\ww_h}\label{acoustic:grad2}
  \end{align}
\end{subequations}
Above, $(\vv_h, e_h)$ is the solution of the state equation
\eqref{acoustic:strong} and $(\ww_h, h_h)$ the solution of the adjoint
equation~\eqref{acoustic:adjweak}.  Since the state equation
\eqref{eq:acoustic} is satisfied in the dG sense,
\eqref{acoustic:grad} simplifies provided $\rho c\tilde c h_h\in P^h$,
or if a quadrature method is used in which the values of
$\rho c\tilde c h_h$ at the quadrature points coincide with the values of a
function in $P^h$ at these points. The latter is, for instance, always
the case when the same nodes are used for the quadrature and the nodal
basis.
Then,
\eqref{acoustic:grad} simplifies to
\begin{subequations}\label{acoustic:grad_simp}
  \begin{align}\label{acoustic:grad_simp1}
    G \tilde c &= -2\nabla(\rho c \tilde c e_h)\cdot \ww_h,\\
    g &= \frac 12 \rho c^2 \diff{e_h}\diff{h_h} +
    \frac 12 \rho \jump{\vv_h}\jump{\ww_h} + 2 \rho
    c\jump{e_h}\cdot\mean{\ww_h}.\label{acoustic:grad_simp2}
  \end{align}
\end{subequations}
As for the one-dimensional advection problem (see
Remark~\ref{rem:adv3}), the upwind flux in the state equation
becomes a downwind flux in the adjoint equation, and thus an upwind
flux for the adjoint equation when solved backwards in time.
\begin{remark}\label{rem:acoustic2}
As in the advection example, the discrete gradient has boundary
contributions that involve jumps of the dG variables at the element
boundaries (see~\eqref{acoustic:grad2}
and~\eqref{acoustic:grad_simp2}). These jumps are at the order of the
dG approximation error and thus tend to zero as the dG solution
converges to the continuous solution either through mesh refinement or
improvement of the approximation on each element.
\end{remark}

\subsection{Maxwell's equations}\label{subsec:maxwell}
Here we derive expressions for the discrete gradient with respect to the
current density in Maxwell's equations (specifically boundary current
density in our case).  This can be used, for instance, in the determination and
reconstruction of antennas from boundary field
measurements~\cite{Nicaise00} and controlling electromagnetic fields
using currents~\cite{Lagnese89, Yousept12}.
The time-dependent Maxwell's equations in a
homogeneous isotropic dielectric domain $\Omega\subset\mathbb R^3$ is given by:
\begin{subequations}\label{eq:maxwells}
  \begin{alignat}{2}
          \mu\mHH_t &=          -  \Curl\mEE && \qquad \text{on}\ \Omega\times(0,T), \label{max:faraday}  \\
     \epsilon\mEE_t &= \phantom{-} \Curl\mHH && \qquad \text{on}\ \Omega\times(0,T), \label{max:ampere}   \\
         \Div\mHH   &= 0                     && \qquad \text{on}\ \Omega\times(0,T), \label{max:gaussmag} \\
         \Div\mEE   &= 0                     && \qquad \text{on}\ \Omega\times(0,T), \label{max:gaussele}
  \end{alignat}
  where $\mEE$ is the electric field and $\mHH$ is the magnetic field.
  Moreover, $\mu$ is the permeability and $\epsilon$ is the
  permittivity, which can both be discontinuous across element
  interfaces.  The impedance $Z$ and the conductance $Y$ of the
  material are defined as
  $Z=\frac{1}{Y}=\sqrt{\frac{\mu}{\epsilon}}$. Note that we follow a
  standard notation for Maxwell's equation, in which the vectors
  $\mHH$ and $\mEE$ are denoted by bold capital letters. Together with
  equations~\eqref{max:faraday}--\eqref{max:gaussele}, we assume the
  initial conditions
  \begin{align}\label{max:ics}
    \mEE(\xx,0) = \mEE_0(\xx), \:\mHH(\xx,0) = \mHH_0(\xx) \quad
    \text{on}\ \Omega,
  \end{align}
  and boundary conditions
  \begin{align}\label{eq:bcs}
    \nn \times \mHH  = -\mJs  \quad \text{on}\ \Gamma.
  \end{align}
  This classic boundary condition can be converted to equivalent
  inflow characteristic boundary conditions~\cite{TengLinChangEtAl08}.
  Here, $\mJs(\xx,t)$ is
  a spatially (and possibly time-dependent) current density flowing
  tangentially to the boundary.
\end{subequations}
If the initial conditions satisfy the divergence
conditions~\eqref{max:gaussmag} and~\eqref{max:gaussele}, the time
evolved solution will as well~\cite{HesthavenWarburton02}.  Thus, the
divergence conditions can be regarded as a consistency condition on
the initial conditions.  We consider a dG discretization of
Maxwell's equations hat only  involves
equations~\eqref{max:faraday} and~\eqref{max:ampere} explicitly.  The dG solution
then satisfies the divergence conditions up to discretization error.
The strong form dG discretization of equation~\eqref{eq:maxwells} is:
Find $(\mHH_h,\mEE_h) \in \mXY$ satisfying the initial conditions
\eqref{max:ics}, such that
\begin{equation}\label{max:strong}
  \begin{split}
   \int_\Omega (\mu     \mHH_{h,t} + \Curl\mEE_h) \cdot\mGG_h \,d\xx
  +\int_\Omega (\epsilon\mEE_{h,t} - \Curl\mHH_h) \cdot\mFF_h \,d\xx=\qquad\qquad\\
  \sum_e\!\int_{\partial\De}\!\!\!\!\! -\left(\nn^-\!\times\!(\mEE_h^\dagger - \mEE_h^-)\right) \cdot\mGG_h \,d\xx
  + \sum_e\!\int_{\partial\De} \!\!\! \left(\nn^-\times(\mHH_h^\dagger - \mHH_h^-)\right) \cdot\mFF_h \,d\xx
  \end{split}
\end{equation}
for all $(\mGG_h,\mFF_h) \in \mXY$, and for all $t\in (0,T)$. The
upwind numerical flux states are given such that
\begin{subequations}\label{max:flux}
  \begin{align}
   \nn^-\times(\mEE_h^\dagger - \mEE_h^-) &=
     -\frac1{2\mean{Y}}\nn^-\times(Y^+\diff{\mEE_h} + \nn^-\times\diff{\mHH_h}), \\
   \nn^-\times(\mHH_h^\dagger - \mHH_h^-) &=
     -\frac1{2\mean{Z}}\nn^-\times(Z^+\diff{\mHH_h} - \nn^-\times\diff{\mEE_h}).
  \end{align}
\end{subequations}
The boundary conditions~\eqref{eq:bcs} are imposed via the upwind numerical flux by setting
exterior values on the boundary of the domain $\Gamma$ such
that
\begin{subequations}
  \begin{align}
    \mHH_h^+ &= -\mHH_h^-+2\mJs, \\
    \mEE_h^+ &=  \mEE_h^-,
  \end{align}
\end{subequations}
with the continuously extended material parameters $Y^+=Y^-$ and $Z^+=Z^-$.
Using the upwind numerical flux implicitly means that the boundary
conditions are only set on the incoming characteristics.

Next, we compute the discrete adjoint equation and the gradient with
respect to the boundary current density $\mJs$ for an objective
functional of the form
\begin{equation}\label{eq:max-cost}
\tilde{\mathcal{J}}(\mJs,\mEE_h) = \int_0^T\!\!\int_\Omega
j_\Omega(\mEE_h)\,d\xx\,dt + \int_0^T\!\!\int_\Gamma r(\mJs)\,d\xx\,dt,
\end{equation}
with differentiable functions $j_\Omega:\Omega\to\mathbb{R}$ and
$r:\Gamma\to\mathbb{R}$.  We introduce the Lagrangian function, and
derive the adjoint equation by imposing that all variations of the
Lagrangian with respect to $\mHH_h$ and $\mEE_h$ must vanish.  After
integration by parts in time $t$, this results in the following
adjoint equation: Find $(\mGG_h,\mFF_h) \in \mXY$ such that
\begin{equation}\label{max:adjweak}
  \begin{split}
   \int_\Omega \mu \mGG_{h,t}\cdot\tilde\mHH_h
  + \mFF_h\cdot(\Curl\tilde\mHH_h) \,d\xx
  +\int_\Omega \epsilon\mFF_{h,t}\cdot\tilde\mEE_h
  -\mGG_h\cdot(\Curl\tilde\mEE_h) \,d\xx=\qquad\qquad\\
  \sum_e\!\int_{\partial\De}\!\!\!\! -\left(\nn^-\times\mFF_h^\dagger\right) \cdot\tilde\mHH_h \,d\xx
  + \sum_e\!\int_{\partial\De}  \left(\nn^-\times\mGG_h^\dagger\right) \cdot\tilde\mEE_h \,d\xx
  - \int_\Omega j_\Omega'(\mEE_h)\cdot\tilde\mEE_h \,d\xx
  \end{split}
\end{equation}
for all $(\tilde\mHH_h,\tilde\mEE_h) \in \mXY$ and all $t\in(0,T)$,
with the final time conditions
\begin{equation}
\mGG_h(\xx,T)=0, \:\mFF_h(\xx,T)=0\quad \text{on}\ \Omega,
\end{equation}
and the adjoint numerical flux states
\begin{subequations}\label{max:adjflux}
  \begin{align}
    \mFF_h^\dagger &= \frac{\mean{Y \mFF_h}}{\mean{Y}} + \frac{1}{2\mean{Y}} \left(\nn^-\times\diff{\mGG_h}\right),\\
    \mGG_h^\dagger &= \frac{\mean{Z \mGG_h}}{\mean{Z}} - \frac{1}{2\mean{Z}} \left(\nn^-\times\diff{\mFF_h}\right),
  \end{align}
\end{subequations}
with exterior values on the boundary $\Gamma$ given by
    $\mGG_h^+ = -\mGG_h^-$ and $\mFF_h^+ =  \mFF_h^-$.
These exterior states enforce the continuous adjoint boundary condition
\begin{equation*}
    \nn \times \mGG  = 0  \quad \text{on}\ \Gamma.
\end{equation*}
To compare with the numerical flux~\eqref{max:flux} of the state equation, we rewrite the
adjoint numerical flux states as
\begin{subequations}\label{max:adjfluxdiff}
  \begin{align*}
   \nn^-\times(\mFF_h^\dagger - \mFF_h^-) &=
     -\frac1{2\mean{Y}}\nn^-\times(Y^+\diff{\mFF_h} - \nn^-\times\diff{\mGG_h}), \\
   \nn^-\times(\mGG_h^\dagger - \mGG_h^-) &=
     -\frac1{2\mean{Z}}\nn^-\times(Z^+\diff{\mGG_h} + \nn^-\times\diff{\mFF_h}).
  \end{align*}
\end{subequations}
Note that, even with discontinuities in the material parameters,
\eqref{max:adjweak} is the weak form of the dG discretization for a
Maxwell's system solved backwards in time.  As in the acoustic
example, the adjoint numerical flux states~\eqref{max:adjflux} come
from the downwind flux.

Differentiating the Lagrangian with respect to the boundary current
$\mJs$ yields an equation for the derivative in direction $\tilde\mJs$, namely
\begin{subequations}\label{max:grad}
\begin{equation}\label{max:grad1}
{\mathcal J}'(\mJs)(\tilde \mJs) = \int_0^T\!\!\!\int_\Gamma r'(\mJs) \cdot \tilde \mJs \,d\xx +
\int_0^T\!\!\int_\Gamma \mg \cdot \tilde \mJs \,d\xx\,dt
\end{equation}
with
\begin{equation}\label{max:grad_parts}
  \mg = \nn^-\times\mFF_h +
         \frac{1}{Y}\left(\nn^-\times\left(\nn^-\times\mGG_h\right)\right),
\end{equation}
\end{subequations}
where $(\mGG_h, \mFF_h)$ is the solution of the discrete adjoint
equation~\eqref{max:adjweak}.
\begin{remark}\label{rem:acoustics}
Since the boundary force $\mJs$ enters linearly in Maxwell's equation,
the gradient expression~\eqref{max:grad} does not involve
contributions from the element boundaries as for the advection and the
acoustic wave example.
\end{remark}

\subsection{Coupled elastic-acoustic wave equation}
\label{subsec:dgea}
Finally, we present expressions for the derivatives with respect to
the primary and secondary wave speeds in the coupled acoustic-elastic
wave equation. This section generalizes Section~\ref{subsec:acoustic}
to the coupled acoustic-elastic wave equation. We derive
derivative expressions with respect to both wave speeds, and allow for
discontinuous wave speeds across elements. We only present a condensed
form of the derivations, and verify our results for the gradient
numerically in Section~\ref{sec:numerics}.

The coupled linear elastic-acoustic wave equation for isotropic
material written in first-order velocity strain form is given as
\begin{subequations}\label{eq:ae}
  \begin{alignat}{2}
    \EE_t & = \frac12\left(\Grad\vv + \Grad\vv^T \right) && \qquad \text{\ on\ } \Omega\times(0,T) \label{eq:ae1}\\
    \rho \vv_t &
      = \Div\left(\lambda\trace(\EE)\ten{I} + 2\mu\EE\right) + \rho \ff && \qquad \text{\ on\ }\Omega\times(0,T) \label{eq:ae2}
  \end{alignat}
  where $\EE$ is the strain tensor, $\vv$ is the displacement
  velocity, $\ten{I}$ is the  identity tensor,
  $\rho=\rho(\xx)$ is the mass density, $\vec{f}$ is a body force per
  unit mass, and $\lambda=\lambda(\xx)$ and
  $\mu=\mu(\xx)$ are the Lam\'e parameters.  In addition to the
  conditions~\eqref{eq:ae1} and~\eqref{eq:ae2} on the body $\Omega$,
  we assume the initial conditions
  \begin{align}
    \vv(0,\xx) = \vv_0(\xx), \: \EE(0,\xx) = \EE_0(\xx), \quad \text{\ for\ }\xx \in \Omega,
  \end{align}
  and the boundary conditions
  \begin{align}
    \SS(\xx,t)\nn = \vec{t}^\text{bc}(t) \qquad \text{on}\ \Gamma.
  \end{align}
Here, $\vec{t}^\text{bc}$ is the traction on the boundary of the body.
\end{subequations}
The stress tensor $\SS$ is related to the strain through the
constitutive relation (here, $\CC$ is the forth-order constitutive tensor):
\begin{equation}
   \SS = \CC\EE = \lambda \trace(\EE) \bs{I} + 2\mu \EE,
\end{equation}
where $\trace(\cdot)$ is the trace operator.  There are also boundary
conditions at material interfaces.  For an elastic-elastic interface
$\Gamma^\text{ee}$ the boundary conditions are
\begin{subequations}
  \begin{alignat*}{2}
    \vv^+  = \vv^-, \quad \SS^+\nn^- & = \SS^-\nn^-  \qquad \text{on}\ \Gamma^\text{ee}
  \end{alignat*}
\end{subequations}
and for acoustic-elastic and acoustic-acoustic interfaces
$\Gamma^\text{ae}$, the boundary conditions are
\begin{subequations}
  \begin{alignat*}{2}
    \nn\cdot\vv^+   = \nn\cdot\vv^-, \quad  \SS^+\nn^- = \SS^-\nn^-  \qquad \text{on}\ \Gamma^\text{ae}.
  \end{alignat*}
\end{subequations}
The strong form dG discretization of equation~\eqref{eq:ae} is:
Find $(\EE_h,\vv_h) \in \aeXY$ such that
\begin{multline}\label{eq:dgae:forward}
  \int_\Omega     \EE_{h,t}:\CC\HH_h   \,d\xx +
  \int_\Omega \rho\vv_{h,t} \cdot\ww_h \,d\xx -
  \int_\Omega \sym(\nabla\vv_h):\CC\HH_h   \,d\xx \\ -
  \int_\Omega \left(\Div(\CC\EE_h) + \vec{f}\right)\cdot\ww_h \,d\xx =
  \sum_e \int_{\partial \De} \sym\left(\nn^-\otimes\left(\vv_h^\dagger-\vv_h^-\right)\right):\CC\HH_h   \,d\xx \\ +
  \sum_e \int_{\partial \De} \left(\left({(\CC\EE_h)}^\dagger - {(\CC\EE_h)}^-\right)\nn^-\right)\cdot\ww_h \,d\xx
\end{multline}
for all $(\HH_h,\ww_h) \in \aeXY$ where $\sym$ is the mapping to get the
symmetric part of a tensor, i.e.,
$\sym(\ten{A})=\frac12\left(\ten{A}+\ten{A}^T\right)$.  Note that the
constitutive tensor $\CC$ is used in the inner product for the weak
form.
Here, the
upwind states are given such that
\begin{subequations}\label{eq:dgae:forwardflux}
  \begin{align}
      \sym\left(\nn^-\otimes\left(\vv_h^\dagger-\vv_h^-\right)\right) &=
    {-k_0\left(\nn^-\cdot\jump{\CC\EE_h} + \rho^+c_p^+ \jump{\vv_h}\right)}
    \left(\nn^-\otimes\nn^-\right) \nonumber\\
    &\quad + k_1
      \sym\left(
        \nn^-\otimes
        \left(\nn^-\times\left(\nn^-\times\jump{\CC\EE_h}\right)\right)
      \right)
      \nonumber\\
    &\quad + k_1\rho^+c_s^+
      \sym\left(
        \nn^-\otimes
        \left(\nn^-\times\left(\nn^-\times\diff{\vv_h}\right)\right)
      \right),
  \\
      \left({(\CC\EE_h)}^\dagger - {(\CC\EE_h)}^-\right)\nn^- &=
    {-k_0\left(\nn^-\cdot\jump{\CC\EE_h} + \rho^+c_p^+ \jump{\vv_h}\right)}
      \rho^-c_p^-\nn^- \nonumber\\
    &\quad + k_1
      \rho^-c_s^-\nn^-\times\left(\nn^-\times\jump{\CC\EE_h}\right)
      \nonumber\\
    &\quad + k_1\rho^+c_s^+
      \rho^-c_s^-\nn^-\times\left(\nn^-\times\diff{\vv_h}\right),
  \end{align}
with $k_0 = 1/(\rho^-c_p^-+\rho^+c_p^+)$ and
\begin{align*}
  k_1 &=
  \begin{dcases}
    \frac{1}{\rho^-c_s^-+\rho^+c_s^+}& \text{when}\ \mu^-\neq0, \\
    0                                & \text{when}\ \mu^-   =0,
  \end{dcases}
\end{align*}
\end{subequations}
where $c_p := \sqrt{(\lambda + 2\mu)/\rho}$ is the primary wave speed
and $c_s := \sqrt{\mu/\rho}$ is the secondary wave speed.  The
traction boundary conditions are imposed through the upwind numerical
flux by setting exterior values on $\Gamma$ to
\begin{align*}
  \vv_h^+ &= \vv_h^-, \\
  \CC^+\EE_h^+\nn^+ &=
\begin{cases}
  -2\vec{t}^\text{bc} + \CC^-\EE_h^-\nn^- &\ \text{if $\mu^-\not=0$},\\
  -2\left(\nn^-\cdot\left(\vec{t}^\text{bc}-\CC^-\EE_h^-\nn^-\right)\right)\nn^-
  - \CC^-\EE_h^-\nn^- &\ \text{if $\mu^+=0$},
\end{cases}
\end{align*}
with the continuously extended material parameters $\rho^+ = \rho^-$,
$\mu^+ = \mu^-$, and $\lambda^+ = \lambda^-$.

We assume a cost function that depends on the primary and secondary
wave speeds $c_p$ and $c_s$ through the solution $\vv_h$ of \eqref{eq:dgae:forward}
\begin{equation}\label{eq:ae-cost}
  {\mathcal{J}}(c_p, c_s) = \int_0^T\!\!\! \int_\Omega
  j_\Omega(\vv_h)\,d\xx\,dt + \int_\Omega r_p(c_p)\,d\xx + \int_\Omega r_s(c_s)\,d\xx.
\end{equation}
By using a sum of spatial Dirac delta distributions in $j_\Omega(\cdot)$,
this can include seismogram data, as common in seismic
inversion.
Using the Lagrangian function and integration by parts in time, we obtain
the following adjoint equation:
Find $(\HH_h,\ww_h) \in \aeXY$ such that
\begin{multline}\label{eq:dgae:adjweak}
  \int_\Omega -\HH_{h,t}:\CC\tEE_h               \,d\xx -
  \int_\Omega \rho\ww_{h,t} \cdot \tvv_h         \,d\xx -
  \int_\Omega \CC\HH_h:\sym(\nabla\tvv_h)            \,d\xx \\ -
  \int_\Omega \ww_h\cdot\left(\Div(\CC\tEE_h)\right) \,d\xx  = -
  \sum_e \int_{\partial \De} \sym(\nn^-\otimes\ww_h^\dagger):\CC\tEE_h\,d\xx \\ -
  \sum_e \int_{\partial \De} \left({(\CC\HH_h)}^\dagger\nn^-\right)\cdot\tvv_h\,d\xx, -
  \int_\Omega j_\Omega(\vv_h)\cdot\tvv_h
\end{multline}
for all $(\tEE_h,\tvv_h) \in \aeXY$ in  with final conditions $\ww_h(T)=0$ and $\HH_h(T)=0$ the fluxes are given by
\begin{align*}
  \sym\left(\nn^-\otimes\ww_h^\dagger\right)
  & =       k_0 \left(\nn^-\cdot\jump{\CC\HH_h} + \nn^-\cdot\left(2\mean{\rho c_p \ww_h}\right) \right)\left(\nn^-\otimes\nn^-\right)\\
  & \quad - k_1\sym\left(\nn^-\otimes\left(\nn^-\times\left(\nn^-\times \jump{\CC\HH_h}\right)\right)\right)\\
  & \quad - k_1\sym\left(\nn^-\otimes\left(\nn^-\times\left(\nn^-\times\left(2 \mean{\rho c_s\ww_h}\right)\right)\right)\right),\\
  {\left(\CC\HH_h\right)}^\dagger\nn^-
  & =
  k_0\left(\nn^-\cdot\left(\left(
      \rho^+c_p^+\CC^-\HH_h^-
    + \rho^-c_p^-\CC^+\HH_h^+\right)\nn^-\right)
    + \rho^-c_p^-\rho^+c_p^+ \jump{\ww_h}\right)\nn^- \\
  & \quad - k_1\nn^-\times\left(\nn^-\times
    \left(\left(\rho^+c_s^+\CC^-\HH_h^- + \rho^-c_s^-\CC^+\HH_h^+\right)\nn^-\right)\right)\\
  & \quad - k_1\rho^-c_s^-\rho^+c_s^+\nn^-\times\left(\nn^-\times \diff{\ww_h}\right).
\end{align*}
We can rewrite this into a form similar to the upwind states of the
state equation~\eqref{eq:dgae:forwardflux} as
\begin{align*}
  \sym\left(\nn^-\otimes\left(\ww_h^\dagger - \ww_h^-\right)\right)
  & =       k_0 \left(\nn^-\cdot\jump{\CC\HH_h} - \rho^+c_p^+\jump{\ww_h}\right)\left(\nn^-\otimes\nn^-\right)\\
  & \quad - k_1\sym\left(\nn^-\otimes\left(\nn^-\times\left(\nn^-\times \jump{\CC\HH_h}\right)\right)\right)\\
  & \quad + k_1\rho^+c_s^+\sym\left(\nn^-\otimes\left(\nn^-\times\left(\nn^-\times \diff{\ww_h}\right)\right)\right),\\
  \left({(\CC\HH_h)}^\dagger - {(\CC\HH_h)}^-\right)\nn^-
  & =       k_0 \left(-\nn^-\cdot\jump{\CC\HH_h} + \rho^+c_p^+\jump{\ww_h}\right)\rho^-c_p^-\nn^-\\
  & \quad + k_1\rho^-c_s^-\nn^-\times(\nn^-\times \jump{\CC\HH_h})\\
  & \quad - k_1\rho^-c_s^-\rho^+c_s^+\nn^-\times(\nn^-\times \diff{\ww_h}).
\end{align*}
Here, the adjoint boundary conditions are
imposed through the adjoint numerical
flux by setting exterior values on $\Gamma$ to
\begin{align*}
  \ww_h^+ &= \ww_h^-, \\
  \CC^+\HH_h^+\nn^+ &=
\begin{cases}
  \phantom{-}\CC^-\HH_h^-\nn^- &\ \text{if $\mu^-\not=0$},\\
           - \CC^-\HH_h^-\nn^- + 2\left(\nn^-\cdot\left(\CC^-\HH_h^-\nn^-\right)\right)\nn^- &\ \text{if $\mu^-=0$},
\end{cases}
\end{align*}
with the continuously extended material parameters $\rho^+ = \rho^-$,
$\mu^+ = \mu^-$, and $\lambda^+ = \lambda^-$.

We assume a discretization of $c_p$ and $c_s$ and a numerical
quadrature rule such that the state equation can be used to simplify
the expression for the gradient; see the discussion in
Section~\ref{subsec:acoustic}.  The discrete gradient
with respect to $c_p$ is then
\begin{subequations}\label{eq:dgae:cpgrad}
\begin{equation}\label{acoustic:cpgrad_gen}
{\mathcal J}_{c_p}(c_p,c_s)(\tilde{c}_p) = \int_\Omega r_p'(c_p)\tilde{c}_p \,d\xx +
\int_0^T\!\!\int_\Omega  G_p\tilde{c}_p \,d\xx\,dt + \sum_{\partial \De} \int_0^T\!\!\int_{\partial \De} g_p\tilde{c}^-_p\,d\xx\,dt,
\end{equation}
where $G_p$ is defined on each element $\De$, and $g_p$ for each
element boundary as follows:
\begin{equation}\label{eq:dgae:cpgrad2}
\begin{split}
  G_p\tilde{c}_p & =
   -2 \left(\nabla\left(\rho\cp\tcp\trace(\EE_h)\right)\right)\cdot\ww_h,
  \\
  g_p &=
           - k_0^2 \rho^- \nn^-\cdot\jump{\CC\EE_h}\nn^-\cdot\jump{\CC\HH_h}
           + k_0^2 \rho^- {\left(\rho^+\cp^+\right)}^2\jump{\vv_h}\jump{\ww_h} \\
    &\quad + k_0^2 \rho^-\rho^+\cp^+\left(\nn^-\cdot\jump{\CC\EE_h}\jump{\ww_h}
           - \jump{\vv_h}\nn^-\cdot\jump{\CC\HH_h}\right) \\
    &\quad + 2 k_0 \rho^- \cp^-
             \trace(\EE_h^-)\left(\nn^-\cdot\jump{\CC\HH_h}-\rho^+\cp^+\jump{\ww_h}\right)\\
    &\quad + 2\rho^-c_p^-\trace(\EE_h^-)\ww_h^-\cdot\nn^-,
\end{split}
\end{equation}
\end{subequations}
where $(\vv_h, \EE_h)$ is the solution of the state equation
\eqref{eq:dgae:forward} and $(\ww_h, \HH_h)$ is the solution of the adjoint
equation~\eqref{eq:dgae:adjweak}.
The discrete gradient of $\mathcal J$ with respect to $c_s$ is
\begin{subequations}\label{eq:dgae:csgrad}
\begin{equation}\label{acoustic:csgrad_gen}
{\mathcal J}_{c_s}(c_p,c_s)(\tilde{c}_s) = \int_\Omega r_s'(c_s)\tilde{c}_s \,d\xx +
\int_0^T\!\!\int_\Omega  G_s\tilde{c}_s \,d\xx\,dt + \sum_{\partial \De} \int_0^T\!\!\int_{\partial \De} g_s\tilde{c}^-_s\,d\xx\,dt,
\end{equation}
where $G_s$ is defined on each element $\De$, and $g_s$ for each
element boundary as follows:
\begin{equation}\label{eq:dgae:csgrad2}
\begin{split}
  G_s\tilde{c}_s & =
   -4 \left(\nabla\cdot\left(\rho c_s\tilde{c}_s\left(\EE_h-\trace(\EE_h)\ten{I}\right)\right)\right)\cdot\ww_h,
  \\
  g_s & = -k_1^2\rho^-\! \left(\nxnx{\jump{\CC\EE_h}}\right) \cdot \left(\nxnx{\big(\jump{\CC\HH_h}-\rho^+c_s^+\diff{\ww_h}\big)}\right) \\
   &\quad - k_1^2\rho^-\rho^+\cs^+\! \left(\nxnx{\diff{\vv_h}}\right) \cdot \left(\nxnx{\big(\jump{\CC\HH_h}-\rho^+c_s^+\diff{\ww_h}\big)}\right) \\
   &\quad + 4k_0\rho^-\cs^- \!\left(\nn^-\cdot\EE_h^-\nn^--\trace(\EE_h^-)\right) \left(\nn^-\cdot\left(\jump{\CC\HH_h}-\rho^+\cp^+\diff{\ww_h}\right)\right) \\
   &\quad + 4k_1\rho^-\cs^-\! \left(\nxnx{\big(\EE_h^-\nn^-\big)}\right)\!\cdot\! \left(\nxnx{\big(\jump{\CC\HH_h}\!-\!\rho^+\cs^+\diff{\ww_h}\big)}\right)\\
   &\quad + 4\rho^-c_s^-\!\left(\left(\EE_h^- - \trace(\EE_h^-)\ten{I} \right)\nn^-\right)\cdot\ww_h^-,
\end{split}
\end{equation}
\end{subequations}
where $(\vv_h, \EE_h)$ is the solution of the state equation
\eqref{eq:dgae:forward} and $(\ww_h, \HH_h)$ is the solution of the
adjoint equation~\eqref{eq:dgae:adjweak}. Note that since we allow for
discontinuous wave speeds $c_p$ and $c_s$ (and perturbations
$\tilde c_p$ and $\tilde c_s$), the boundary contributions to the gradients,
i.e., the last terms in~\eqref{acoustic:cpgrad_gen} and
\eqref{acoustic:csgrad_gen} are written as sums of integrals over
individual element boundaries, i.e., each boundary face appears twice
in the overall sum. This differs from the previous examples, where we
assumed continuous parameters and thus combined contributions from
adjacent elements to shared faces $\mathsf f$.
\begin{remark}\label{rem8:dgea}
The expressions for the $c_p$-gradient~\eqref{eq:dgae:cpgrad} reduce
to the result found for the acoustic
equation~\eqref{acoustic:grad_gen} and~\eqref{acoustic:grad_simp} for
continuous parameter fields $\rho, c_p, c_s$ and continuous parameter
perturbations $\tilde c_p$. To verify this, one adds contributions
from adjacent elements to common boundaries, and the terms
in~\eqref{eq:dgae:cpgrad2} combine or cancel.

\end{remark}
\begin{remark}\label{rem7:dgea}
Above, we have derived expressions for the derivatives with respect to
the primary and secondary wave speeds. If, instead of $c_p$ and $c_s$,
derivatives with respect to an alternative pair of parameters in the
stress tensor---such as the Lam\'e parameters, or Poisson's ratio and
Young's modulus---are derived, the adjoint equations remain
unchanged, but the expressions for the derivatives change according
to the chain rule.
\end{remark}

\section{Numerical verification for coupled elastic-acoustic wave propagation}\label{sec:numerics}

\begin{figure}[ht]\centering
  \begin{minipage}{.67\columnwidth}
  \begin{tikzpicture}[scale=.95]
    \begin{axis}[xlabel=mesh level, ylabel=directional derivative, only marks,
        xtick={1,2,3}, xmin=0.5, xmax=3.5]
      \addplot [color=yellow!90!black, mark=*, mark size=3pt] coordinates {%
        (1, 1.2443581984E+00)
        (2, 9.6318263350E-01)
        (3, 9.5582715331E-01)
      };
      \addlegendentry{$d^\text{fd}_\epsilon, \epsilon=10^{-3}$}
      \addplot [color=red!40!yellow, mark=*,  mark size=3pt] coordinates {%
        (1, 1.2198414986E+00)
        (2, 9.3174411305E-01)
        (3, 9.2088677889E-01)
      };
      \addlegendentry{$d^\text{fd}_\epsilon, \epsilon=10^{-4}$}
      \addplot [color=red!80!yellow, mark=*,  mark size=3pt] coordinates {%
        (1, 1.2173892827E+00)
        (2, 9.2860042571E-01)
        (3, 9.1739423311E-01)
      };
      \addlegendentry{$d^\text{fd}_\epsilon, \epsilon=10^{-5}$}
      \addplot [color=blue, mark=square,  mark size=2.5pt, line width=0.7pt] coordinates {%
        (1 , 1.2171168922E+00)
        (2, 9.2825094400E-01)
        (3, 9.1700439384E-01)
      };
      \addlegendentry{disc.~grad.~$d^\text{di}$}
      \addplot [color=black, mark=diamond,  mark size=3.5pt, line width=0.7pt] coordinates {%
        (1 , 1.2945884839E+00)
        (2, 8.7662282604E-01)
        (3, 9.1127305367E-01)
      };
      \addlegendentry{cont.~grad.~$d^\text{co}$}
    \end{axis}
  \end{tikzpicture}
  \end{minipage}
  \hfill
  \begin{minipage}{.32\columnwidth}
  \begin{tabular}[b]{|c|c|}
    \hline
    \multicolumn{2}{|c|}{mesh level 1} \\ \hline
     $d^\text{fd}_{\epsilon},  \epsilon=10^{-2}$ & {\bf 1}.489022\\
     $d^\text{fd}_{\epsilon},  \epsilon=10^{-3}$ & {\bf 1.2}44358\\
     $d^\text{fd}_{\epsilon},  \epsilon=10^{-4}$ & {\bf 1.21}9841\\
     $d^\text{fd}_{\epsilon},  \epsilon=10^{-5}$ & {\bf 1.217}389\\
     $d^\text{fd}_{\epsilon},  \epsilon=10^{-6}$ & {\bf 1.2171}43\\
     $d^\text{fd}_{\epsilon},  \epsilon=10^{-7}$ & {\bf 1.21711}8\\ \hline
     $d^\text{di}                            $ & 1.217117\\ \hline
  \end{tabular}
  \vspace{6mm}
  \end{minipage}
  \caption{Directional derivatives computed using one-sided finite
    differences \eqref{eq:dgae:fd}, and the discrete and the
    continuous gradients \eqref{eq:dgae:grad}. Left: Results on mesh
    levels 1,2 and 3 corresponding to meshes with 16, 128 and 1024
    finite elements with polynomial order $N=4$. The finite difference
    directional derivatives $d^\text{fd}_\epsilon$ converge to the
    discrete gradient $d^\text{di}$ as $\epsilon$ is reduced. Note
    that as the mesh level is increased, the continuous gradient
    $d^\text{co}$ converges to $d^\text{di}$. Right: Convergence of
    finite difference directional derivative on the coarsest
    mesh. Digits for which the finite difference gradient coincides
    with the discrete gradient are shown in bold.\label{fig:dgae}}
\end{figure}
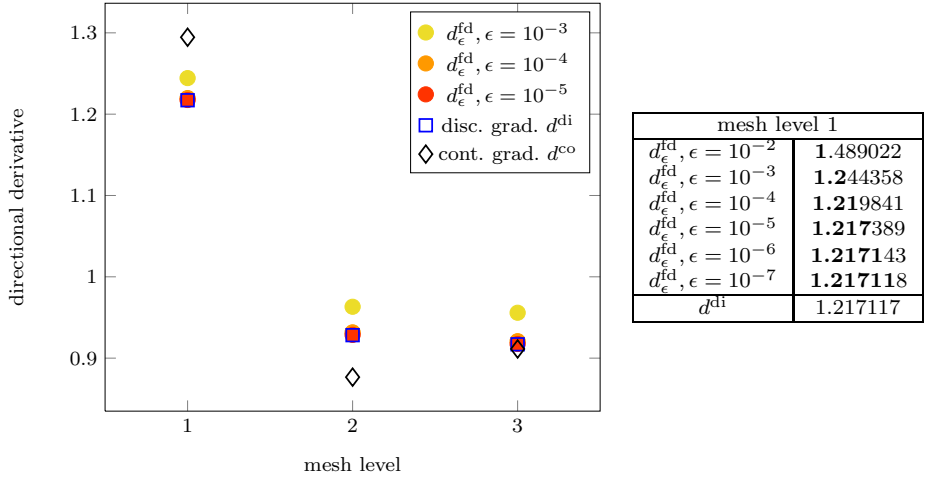

Here, we numerically verify the expressions for the discrete gradients
with respect to the wave speeds for the elastic-acoustic wave problem
derived in Section~\ref{subsec:dgea}. For this purpose, we compare
directional finite differences with directional gradients based on the
discrete gradient. To emphasize the correctness of the discrete
gradient, we use coarse meshes in these comparisons, which
underresolve the wave fields.  As test problem, we use the Snell law
example from Section~6.2 in~\cite{WilcoxStadlerBursteddeEtAl10} with
the material parameters and the wave incident angle specified
there. For our tests, we use the simple distributed objective function
\begin{equation*}
\mathcal J(c_p,c_s) := \int_0^T\!\!\int_\Omega \vv_h\cdot\vv_h\,d\xx\,dt
\end{equation*}
The discretization of the wave equation follows
\cite{WilcoxStadlerBursteddeEtAl10}, i.e., we use spectral elements
based on Gauss-Lobatto-Lagrange (GLL) points on hexahedral meshes. The
use of GLL quadrature results in underintegration even if the elements
are images of the reference element under an affine transformation.
In Figure~\ref{fig:dgae}, we
summarize results for the directional derivatives in the direction
$\tilde c_p:=\sin(\pi x)\cos(\pi y)\cos(\pi z)$; we compare the finite
difference directional derivatives
\begin{equation}\label{eq:dgae:fd}
d_\epsilon^\text{fd}:=\frac{\mathcal J(c_p+\epsilon\tilde c_p) -
  \mathcal J(c_p)}{\epsilon}
\end{equation}
with the directional derivatives $d^\text{di}$ and $d^\text{co}$
defined by
\begin{equation}\label{eq:dgae:grad}
d^\text{di} :=\mathcal J_{c_p}(c_p,c_s)(\tilde c_p),\quad
d^\text{co} :=\mathcal J^\text{cont}_{c_p}(c_p,c_s)(\tilde c_p),
\end{equation}
where $\mathcal J_{c_p}(c_p,c_s)$ denotes the discrete gradient
\eqref{eq:dgae:cpgrad}, and $\mathcal J^\text{co}_{c_p}(c_p,c_s)$
denotes the gradient obtained when neglecting the jump term in the
boundary contributions $g_p$ in~\eqref{eq:dgae:cpgrad2}. These jump
terms are likely to be neglected if the continuous gradient
expressions are discretized instead of following a fully discrete
approach. The resulting error is of the order of the discretization
error and thus vanishes as the discrete solutions converge. However,
this error can be significant on
coarse meshes, on which the wave solution is not well resolved.

Next, we study the accuracy of the discrete gradient for pointwise
perturbations to the wave speed. Since the same discontinuous basis
functions as for the wave equation are also used for the local wave
speeds, a point perturbation in $c_p^-$ or $c_s^-$ at an element
boundary face results in a globally discontinuous perturbation
direction $\tilde c_p$ and $\tilde c_s$.  In Table~\ref{tab:compare},
we present the discrete directional gradient $d^\text{di}$ with finite
difference directional gradients $d^\text{fd}_\epsilon$ for unit
vector perturbations of both wave speeds.  Compared to in the table in
Figure~\ref{fig:dgae}, where the directional derivatives for smooth
perturbations are reported, pointwise perturbations of the wave speeds
$c_p$ or $c_s$ result in smaller changes in the cost functional, and
numerical roundoff influences the accuracy of finite difference
directional derivatives. As a consequence, fewer digits coincide
between the finite difference directional derivatives and the discrete
gradients.

\begin{table}\centering
  \caption{Comparison of pointwise material gradients for Snell
    problem from~\cite[Section 6.2]{WilcoxStadlerBursteddeEtAl10}.
    The derivatives $d^\text{fd}_\epsilon$ and $d^\text{di}$ with
    respect to the local wave speed (either $c_p$ or $c_s$) for points
    with coordinates $(x,y,z)$ are reported. We use the final time
    $T=1$ and spectral elements of polynomial order $N=6$ in
    space. The meshes for level 1 and 2 consist of 16 and 128 finite
    elements, respectively.  Digits where the finite difference
    approximation coincides with the discrete gradient are shown in
    bold.\label{tab:compare}}
  \begin{tabular}{|c|c|c|c|ccc|}
    \hline
 mesh level & $(x,y,z)$ & \!\!pert.\!\! & \!\!$d^\text{di}$\!\! & \multicolumn{3}{c|}{$d^\text{fd}_\epsilon$} \\% & $d^\text{co}$ \\
$\#$tsteps   &           & \!\!field\!\! &                 & \!\!$\epsilon=10^{-3}$\!\! & \!\!$\epsilon=10^{-4}$\!\! & \!\!$\epsilon=10^{-5}$\!\! \\  \hline
 1/101        & $(0,0,0)$       & $c_p$ & 1.8590e-4 &  {\bf 1.85}49e-4 & {\bf 1.85}81e-4& {\bf 1.85}60e-4 \\%& 1.1054e-04 \\
 2/202        & $(0,0,0)$       & $c_p$ & 2.2102e-5 &  {\bf 2.2}094e-5 & {\bf 2}.2007e-5& {\bf 2}.1504e-5  \\%& 2.3655e-5  \\
 1/101        & $(0,0,1)$       & $c_s$ & 1.1472e-5 &  {\bf 1.14}53e-5 & {\bf 1.1}372e-5 & {\bf 1}.0942e-5   \\
 1/101        & $(-0.5,-0.5,0.5)$ & $c_s$ & 2.8886e-3 & {\bf 2.88}02e-3 & {\bf 2.88}77e-3 & {\bf 2.88}70e-3   \\
    \hline
  \end{tabular}
\end{table}

\section{Conclusions}\label{sec:conclusions}
Our study yields that the discretely exact adjoint PDE of a
dG-discretized linear hyperbolic equation is a proper dG
discretization of the continuous adjoint equation, provided an upwind
flux is used. Thus, the adjoint PDE converges at the same rate as the
state equation. When integration by parts is avoided to eliminate
quadrature errors, a weak dG discretization of the state PDE leads to
a strong dG discretization of the adjoint PDE, and vice versa. The
expressions for the discretely exact gradient can contain
contributions at element faces and, hence, differ from a straightforward
discretization of the continuous gradient expression.  These element
face contributions are at the order of the discretization order and are
thus more significant for poorly resolved state PDEs. We believe that
these observations are relevant for inverse problems and optimal
control problems governed by hyperbolic PDEs discretized by the
discontinuous Galerkin method.

\section*{Acknowledgments}
 We would like to thank Jeremy Kozdon and Gregor Gassner for fruitful
 discussions and helpful comments, and Carsten Burstedde for his help
 with the implementation of the numerical example presented in
 Section~\ref{sec:numerics}.  Support for this work was provided
 through the U.S.~National Science Foundation (NSF) grant
 CMMI-1028889, the Air Force Office of Scientific Research's
 Computational Mathematics program under the grant FA9550-12-1-0484,
 and through the Mathematical Multifaceted Integrated Capability
 Centers (MMICCs) effort within the Applied Mathematics activity of
 the U.S.~Department of Energy's Advanced Scientific Computing
 Research program, under Award Number DE-SC0009286.  The views
 expressed in this document are those of the authors and do not reflect
 the official policy or position of the Department of Defense or the
 U.S.~Government.
 \bibliography{ccgo}

\end{document}